\documentclass [12pt]{amsart}
\usepackage{latexsym,amssymb,amsmath}

\usepackage{simplemargins}
\setleftmargin{2.75cm}
\setrightmargin{2.75cm}
\settopmargin{2.0in}
\setbottommargin{2.0in}
\setallmargins{2.75cm}

\def\pf{\noindent\emph{Proof: }}
\def\stop{\hfill$\Box$}

\usepackage{amsmath, amsfonts, amssymb, amsthm}
\newtheorem{thm}{Theorem}

\newtheorem{lemma}[thm]{Lemma}

\newtheorem{defn}[thm]{Definition}
\newtheorem{prop}[thm]{Proposition}

\numberwithin{thm}{section}

\DeclareMathOperator{\Vol}{Vol}

\DeclareMathOperator{\Ric}{Ric}
\DeclareMathOperator{\Hess}{Hess}
\DeclareMathOperator{\diver}{div}

\begin{document}

\title [ Scalar Curvature Deformation Equation ] {The Scalar Curvature Deformation Equation on Locally Conformally Flat Manifolds }

\author{Yu Yan}
\address{Department of Mathematics\\
         The University of British Columbia\\
         Vancouver, B.C., V6T 1Z2\\ Canada}
   \email{yyan@math.ubc.ca}

\begin{abstract}We study the equation $\Delta _g u -\frac{n-2}{4(n-1)}R(g)u+Ku^p=0 \,\ (1+\zeta
\leq p \leq \frac{n+2}{n-2})$ on locally
conformally flat compact manifolds $(M^n,g)$.  We prove the following: (i) When the scalar curvature $R(g)>0$ and the dimension $n \geq 4$, under suitable conditions on $K$, all positive solutions $u$ have uniform upper and lower bounds; (ii) When the scalar curvature $R(g)\equiv 0$ and $n \geq 5$, under suitable conditions on $K$, all positive solutions $u$ with bounded energy have uniform upper and lower bounds.  We also give an example to show that the energy bound condition for the uniform estimates in \cite{YY3} is necessary.

\end{abstract}

\maketitle

 \section{Introduction}

Let $(M^n,g)$ be an n-dimensional compact manifold with metric $g$,
and we use $R(g)$ to denote the scalar curvature of $g$.  Let $u$ be a positive function defined on $M$.  The scalar curvature of the conformally deformed metric $u^{\frac{4}{n-2}}g$ is given by
$$R(u^{\frac{4}{n-2}}g)=-c(n)^{-1}u^{-\frac{n+2}{n-2}} \big (\Delta
_{g} u-c(n)R(g)u \big ) \hspace{.15in} \text{where }
c(n)=\frac{n-2}{4(n-1)}.$$

\noindent
The Yamabe Theorem, which was proved by the work of Trudinger \cite{Tr}, Aubin \cite{Au} and
Schoen \cite{S1}, says that there exists $u>0$ such that $R(u^{\frac{4}{n-2}}g)$ is equal to some constant $K$.  The P.D.E. formulation of this theorem is that the equation

$$
\Delta_{g}u - c(n) R(g)u +c(n)Ku^{\frac{n+2}{n-2}}=0
$$
has a positive solution for some constant $K$.

\noindent
In \cite{ES}, J. Escobar and R. Schoen extended this result to the case when $K$ is a function on $M$.  They proved that under certain conditions on $K$, the above equation has a positive solution $u$ when $R(g)>0$ or $R(g) \equiv 0$.

In fact, in those existence results the solution minimizes the associated constraint variational problem and can be obtained as a limit of a sequence of solutions of the corresponding subcritical equations.  Therefore, a natural question is whether non-minimal solutions can also be produced from solutions of the subcritical equations.  We would like to know if there are uniform estimates for solutions of the equation

\begin{equation}
\label{eq:yamabe}
\Delta_{g}u - c(n) R(g)u +Ku^p=0  \hspace{.3in} \text{ where } \hspace{.1in} 1+\zeta \leq p \leq \frac{n+2}{n-2}.
\end{equation}

\noindent
This was proved to be true by R. Schoen \cite{S2, S6} when $K$ is a positive constant, $R(g)>0$, and $(M^n,g)$ is locally conformally flat and not conformally diffeomorphic to $S^n$.  By the work of Y. Li and M. Zhu \cite{LZ}, this is also true when $K$ is a positive function on a 3-dimensional compact manifold $(M^3,g)$ which has $R(g)>0$ and is not conformally diffeomorphic to $S^3$.  In the case when $K$ is a positive constant, this result by Li and Zhu was extended to dimensions $n=4,5$ by O. Druet in \cite{Druet1, Druet2}.  Then it was extended further to dimensions $n \leq 7$ independently by Y. Li and L. Zhang \cite{Li-Zhang} and F.C. Marques \cite{Marques}; when the dimension $n \geq 8$, it was proved to be true by Li and Zhang \cite{Li-Zhang} under an additional assumption on the Weyl tensor of the backgroud metric $g$.

 In \cite{YY3} we proved uniform estimates for solutions with bounded energy when $K$ is a function satisfying certain conditions on a 3 or 4 dimensional locally conformally flat manifold with zero scalar curvature.  In this paper we study this problem when $K$ is a function on locally conformally flat manifolds $(M^n,g)$.  We consider two separate cases: $R(g)>0$ and $R(g) \equiv 0$.

\subsection{Manifolds of Zero Scalar Curvature}
\label{subsection:scalar-flat-result}

When the scalar curvature $R(g)\equiv 0$ on the manifold $M$, equation  (\ref {eq:yamabe}) becomes
\begin{equation}
\label{eq:main}
\Delta _g u + Ku^p=0 
\hspace{.3in} \text{ where } 1+\zeta \leq p \leq \frac{n+2}{n-2}.
\end{equation}

\noindent
The necessary conditions for the existence of a solution $u>0$ are that $K$ changes sign on $M$ and $\int _{M} K dv_g <0$.

 The corresponding existence result is the following theorem in \cite{ES}:

\begin{thm}
\label {thm:ESflat}
{\rm (Escobar--Schoen \cite{ES}).}
Suppose $M$ is locally conformally flat with zero scalar
curvature. Suppose $K$ is a nonzero smooth function on $M$ satisfying
the condition that there is a maximum point $P_0 \in M$ of $K$ at which all
derivatives of $K$ of order less than or equal to $(n-3)$ vanish. Then
$K$ is the scalar curvature of a metric $\bar{g}=u^{\frac{4}{n-2}}g$
for some $u>0$ on $M$ if and only if $K$ satisfies\\
(i)  $K$ changes sign\\
(ii) $\int _{M} K dv_g <0$.\\
When the dimension $n=3,4$, the flatness condition on $K$ is automatically satisfied and the
locally conformally flat assumption on $M$ can be removed.
\end{thm}

In \cite{YY3}, we proved a compactness theorem when the dimension of $M$ is equal to 3 or 4.
\begin{thm}
\label{thm:lowdim}
{\rm (\cite{YY3}).}
Let $(M,g)$ be a three or four dimensional locally conformally flat compact
manifold with $R(g) \equiv 0$. Let $\mathcal{K}:=\{ K \in C^3(M): K>0 $
somewhere on $M, \int_M K dv_g\leq -{C_K}^{-1}<0,$ and $  \| K
\|_{C^3(M)}\leq C_K \}$ for some constant $C_K$, and
$S_{\Lambda} :=\{u: u>0 $ solves $(\ref{eq:main})$ with $K \in \mathcal{K}, $ and $ E(u):=\int _M |\nabla u|^2 dv_g
\leq \Lambda \}.$ Then there exists $C=C(M, g, C_K, \Lambda, \zeta)>0$ such
that $u \in S_{\Lambda}$ satisfies $ \|u\|_{C^3(M)}\leq C$ and $\displaystyle \min _M u \geq C^{-1}$.
\end{thm}

\noindent
In Section \ref {section:example-notations} we will give an example which shows that these estimates {\it cannot } be improved to be independent of the energy $E(u)$.

Next we give a similar theorem on manifolds of dimension $n \geq 5$.  We first need to define a flatness condition on $K$ as follows.

\begin{defn}
\label{defn:flatonK}
A function $K \in C^{n-2}(M)$ is said to satisfy the {\bf flatness condition $(*)$} if near each critical point $P$ of $K$ where $K(P)>0$, there exist a neighborhood and a constant $C_0$ such that in that neighborhood
$$|\nabla ^p K| \leq C_0  | \nabla K |^{\frac{n-2-p}{n-3}} \hspace{.2in} \text{ for } \hspace{.2in}   2\leq p \leq n-3,$$
where $\nabla ^p K $ is the $p$-th covariant derivative of $K$.
\end{defn}

\noindent
Note that this implies in particular all partial derivatives of $K$ up to order $n-3$ vanish at those critical points, and the order of flatness is the same as that in Theorem \ref {thm:ESflat}.  A simple example of a function satisfying this condition is a function which can be expressed near the critical points as $K(z)=a+b|z|^{n-2}$, where $a, b$ are two constants and $z$ is a local coordinate system centered at the critical point.  This type of flatness condition also appeared in \cite{YLi1} and \cite{YLi2}, where Y. Li studied the problem of prescribing scalar curvature functions on $S^n$.

We are ready to state the theorem: 

\begin{thm}
\label {thm:main}
Let $(M^n,g)$ be a locally conformally flat compact manifold with $R(g)\equiv 0$, and its dimension $n\geq 5$. Let $K \in C^{n-2}(M)$ be a function which satisfies the flatness condition $(*)$; additionally, $K$ is positive somewhere on $M$ and $\int_M K dv_g <0$. If $u$ is a positive solution of equation $(\ref{eq:main})$ with bounded energy $E(u):=\int_M |\nabla u|^2 dv_g \leq \Lambda$, then there exists a positive constant $C$ such that $ \|u\|_{C^3(M)}\leq C$ and $\displaystyle \min _M u \geq C^{-1}$, where $C$ depends on $ M, g, \|K\|_{C^{n-2}(M)}, \int_M K dv_g, \Lambda,$ and $ \zeta$.
\end{thm}

\subsection{Manifolds of Positive Scalar Curvature}
\label {subsection: scalar-pos-result}

When the scalar curvature $R(g)>0$, the necessary condition for equation (\ref {eq:yamabe}) to have a positive solution is that $K>0$ somewhere on the manifold.  The following existence result was proved in \cite {ES}.

\begin{thm}
\label {thm:ESpositive}
{\rm (Escobar--Schoen \cite{ES}).}
Suppose $M$ is a locally conformally flat manifold with positive scalar
curvature which is not simply connected, and $K$ is a smooth
function on $M$ which is somewhere positive, and there is a maximum
point $P_0$ of $K$ at which all partial derivatives of $K$ of order
less than or equal to $(n-2)$ vanish. Then equation $(\ref {eq:yamabe})$ has a
positive solution.\\
When the dimension $n=3$, the flatness condition on $K$ is automatically satisfied and the locally conformally flat assumption on $M$ can be removed.
\end{thm}

The compactness result when $n=3$ was proved in \cite{LZ}.

\begin{thm}
\label {thm:Li-Zhu}
{\rm (Li--Zhu \cite{LZ}).} Let $(M,g)$ be a three dimensional smooth
compact Riemannian manifold with positive scalar curvature which is not conformally equivalent to the
standard $S^3$. Then for any $1<p\leq 5$ and positive function $K \in
C^2(M)$, there exists some constant $C$ depending only on $M,g,
\|K\|_{C^2(M)}$, and the positive lower bound of $K$ and $p-1$ such
that  
$$      \frac{1}{C} \leq u \leq C \hspace{.3in} \text {and} \hspace {.3in}
\|u\|_{C^3(M)}\leq C  $$ 
for all positive solutions u of $\Delta _g u - \frac {1}{8}R(g) u
+Ku^p =0$.
\end{thm}

We will give a compactness theorem when the dimension $n \geq 4$.  But $K$ needs to satisfy a flatness condition near its critical points. 
\begin{defn}
\label{defn:flatonK-pos}
A function $K \in C^{n-1}(M)$ is said to satisfy the {\bf flatness condition} $(**)$ if near each critical point of $K$, 
there exist a neighborhood and a constant $C_0$ such that in that neighborhood
$$|\nabla ^p K| \leq C_0  | \nabla K |^{\frac{n-1-p}{n-2}} \hspace{.2in} \text{ for } \hspace{.2in}   2\leq p \leq n-2,$$
where $\nabla ^p K $ is the $p$-th covariant derivative of $K$.
\end{defn}

\noindent
Under this condition all partial derivatives of $K$ up to order $n-2$ vanish at the critical points, which is consistent with the condition given in Theorem \ref {thm:ESpositive}. A simple example of a function satisfying this condition is a function which can be expressed near the critical points as $K(z)=a+b|z|^{n-1}$, where $a, b$ are two constants and $z$ is a local coordinate system centered at the critical point.

Our theorem is:

\begin{thm}
\label {thm:main-pos}
Let $(M^n,g)$ be a locally conformally flat compact manifold with $R(g)>0$.  Assume $M$ is not conformally diffeomorphic to $S^n$, and its dimension $n \geq 4$. Let $K \in C^{n-1}(M)$ be a positive function which satisfies the flatness condition $(**)$.  There exists a positive constant $C$ such that $\|u\|_{C^3(M)}\leq C$ and $\displaystyle \min _M u \geq C^{-1}$ for any positive solution $u$ of equation $(\ref {eq:yamabe})$, where $C$ depends on $M, g, \zeta$ and $ \|K\|_{C^{n-1}(M)}$.
\end{thm}

\noindent
Note that because we assume $K>0$ in this theorem, there is no assumption on the energy of $u$, which was introduced in the scalar-flat case to overcome the difficulty caused by the sign changing of $K$. 

\section {The Example and Some Notations}
\label {section:example-notations}

Let $(M^n,g)$ be a compact manifold with $R(g) \equiv 0$ and $n=3$ or $4$. (In fact in this example $M$ does not need to be locally conformally flat.)  We choose $K \in C^3(M)$ satisfying the following conditions:
\begin{itemize}
\item $K>0$ somewhere on $M$,
\item $\int _M K dv_g \leq -C_K^{-1}<0$ and $\| K \| _{C^3(M)} \leq C_K$, where $C_K$ is a positive constant,
\item the set $\{ x \in M : K(x)=0 \} = \overline{U}$ for some open set $U \subset M$. 
\end{itemize}

We define 
\begin{displaymath}
K_i(x) = \left \{ \begin{array}{ll}
\frac{K(x)}{i} &  \text{ if } K(x) >0\\
K(x) & \text{ if } K(x) \leq 0
     \end{array} \right.        
\end{displaymath}

\noindent
Since on $\partial U$ all derivatives of $K$ up to order $3$ are zero, it follows that $K_i \in C^3(M)$.  Furthermore, by this definition $K_i \in \mathcal{K}$, where $\mathcal{K}$ is as defined in Theorem \ref {thm:lowdim}.  Then by Theorem \ref {thm:ESflat} there exists $u_i>0$ which satisfies $\Delta _g u_i + K_i u_i^{\frac{n+2}{n-2}}=0$.

Now suppose there is a constant $C$ independent of $i$ such that $\displaystyle \max _M u_i \leq C$.  As proved in Section 2 of \cite {YY3}, this implies that $\{ u_i \} $ is uniformly bounded away from $0$ and $\| u_i \| _{C^3(M)}$ is bounded above uniformly.  Then passing to a subsequence $\{u_i\}$ converges in the $C^2$-norm to a function $u>0$, and $u$ satisfies $\Delta _g u + \tilde{K} u^{\frac{n+2}{n-2}}=0$ where 
\begin{displaymath}
\tilde{K}(x) \,\, = \,\, \lim _{i \to \infty} K_i(x) \,\, = \,\, \left \{ \begin{array}{ll}
0 &  \text{ if } K(x) >0\\
K(x) & \text{ if } K(x) \leq 0
     \end{array} \right.        
\end{displaymath}

\noindent
However, because $\tilde{K}$ is nowhere positive and somewhere negative, the equation  $\Delta _g u + \tilde{K} u^{\frac{n+2}{n-2}}=0$ cannot have a positive solution by Theorem \ref {thm:ESflat}.  This contradiction shows that estimates like the ones in Theorem \ref {thm:lowdim} can not be true without the energy bound assumption on $u$. 

\vspace{.2in}

Next we prove Theorems \ref {thm:main} and \ref {thm:main-pos}.  We will prove Theorem \ref {thm:main} in Sections \ref {section:prelim-flat} to \ref {section:case2}, and the proof of Theorem \ref {thm:main-pos} will be given in Section \ref {section:positive}.   We first give some definitions and a lemma which will be used in both proofs. 

\begin{defn}
\label{defn:bupt}
We call a point $\bar{x}$ on a manifold $M$ a {\bf blow-up point} of a sequence $\{ u_i \}$ if $\bar{x}=\displaystyle \lim _{i
\rightarrow \infty} x_i$ for some $\{ x_i \} \subset M$ and $u_i(x_i) \rightarrow \infty$.
\end{defn}

\begin{defn}
\label{defn:iso}
Suppose $u_i$ satisfies $\Delta _{g_i} u_i-c(n)R(g_i)u_i+K_iu_i^{p_i}=0$, where $\{g_i\}$ converges to some metric $g_0$.  A point $\bar{x} \in M$ is called an {\bf isolated blow-up
point} of $\{ u_i \}$ corresponding to $\{g_i\}$ if there exist local maximum points $x_i$ of
$u_i$ and a fixed radius $r_0>0$ such that 
\begin{itemize}
\item $x_i \rightarrow \bar{x}$,
\item $u_i(x_i) \rightarrow \infty$,
\item $u_i(x) \leq C \left (d_{g_i} (x, x_i)\right )^{-\frac{2}{p_i-1}}$ \hspace{.05in} for any $x \in B_{r_0}(x_i)$, where the constant $C$ is independent of $i$.
\end{itemize}
\end{defn}

\begin{lemma}
\label{lemma:simpleharnack}
If $\bar{x}=\displaystyle \lim _{i \rightarrow \infty} x_i$ is an isolated blow-up point
of $\{u_i\}$ corresponding to $\{g_i\}$, and $ K_i$ is uniformly bounded, then there exists a constant $C$ independent of $i$ and
$r$ such that $$\max _{\partial B_r(x_i)} u_i(x) \leq C \min _{\partial
B_r(x_i)} u_i(x)$$
for any $0<r \leq r_0$.
\end{lemma}
This can be proved as in \cite {YY3} in the proof of Lemma 5.2.

\begin{defn}
\label{defn:simple}
$\bar{x}$ is called a {\bf simple blow-up point} of $\{u_i\}$ if it is an isolated blow-up point and there exists $\bar{r} >0$ independent of $i$ such that
$\bar{w}_i(r)$ has only one critical point for $r \in (0, \bar{r})$.  Here $\bar{w}_i(r):=r^{\frac{2}{p_i-1}}\bar{u}_i(r)=\Vol (S_r)^{-1} \int _{S_r} |z|^{\frac{2}{p_i-1}}u_i(z)d \Sigma _g$ and $z$ is the conformally flat coordinate system centered at each $x_i$.
\end{defn}

\section{Initial Steps of the Proof of Theorem \ref {thm:main}.}
\label{section:prelim-flat}

The proof of Theorem \ref{thm:main} follows along the same line of reasoning as the proof of Theorem \ref{thm:lowdim}, which is done in \cite{YY3}.  As proved in Section 2 of \cite{YY3}, a lower bound on $u$ follows directly if there is a uniform upper bound on $u$.  By the standard elliptic theory and Sobolev embedding theorem, a bound on the $C^0$-norm of $u$ easily implies a bound on its $C^3$-norm.  Therefore, to prove Theorem \ref {thm:main} we only need to show that there is a uniform upper bound on $u$.

By an argument identical to that in Section 3 of \cite{YY3}, we can show that there exists a positive constant $\eta = \eta (M, g, n, \| K \|_{C^{n-2}(M)}, \Lambda)$, such that on the set $K_{\eta} := \{x \in M : K(x) <\eta\}$, $u$ has a uniform upper bound depending only on $M, g, n, \| K \|_{C^{n-2}(M)},$ and $\Lambda$.  Thus it is left to show that $u$ is uniformly bounded on the set where $K \geq \eta$.  
We have the following proposition.
\begin{prop}
\label{prop:iso}
Given $\epsilon >0, R >> 0$, there exists $C=C(\epsilon, R)$ such that
if $u$ is a solution of equation $(\ref{eq:main})$ and
$$\max _{x \in M} \Big( \big(d_g(x, K_{\frac{\eta}{2}})\big
)^{\frac{2}{p-1}}u(x) \Big) > C,$$
then there exists $\{ x_1, ..., x_N \} \subset M \setminus
K_{\frac{\eta}{2}}$ with $N $ depending on $u$, and

\begin{itemize}

\item Each $x_i$ is a local maximum point of $u$ and the geodesic balls $\{
B_{\frac{R}{u(x_i)^{\frac{p-1}{2}}}} (x_i)\}$ are disjoint. 

\item $|\frac{n+2}{n-2}-p | < \epsilon$ and in the coordinate system $y$ so
chosen that $z=\frac{y}{u(x_i)^{\frac{p-1}{2}}}$ is the conformally flat coordinate system centered at $x_i$, we have
 
$$ \Bigg \| u(x_i)^{-1}u\left (\frac
{y}{u(x_i)^{\frac{p-1}{2}}} \right ) -\bar{v}(y) \Bigg \|_{C^2(B_{2R}(0))} <
\epsilon$$ on the ball $B_{2R}(0) \subset \mathbf{R}^n(y)$, where
$$\bar{v} (y)=\left (1+
\frac{K(x_i)}{n(n-2)}|y|^2 \right )^{-\frac{n-2}{2}}.$$ 

\item There exists $C=C(\epsilon, R)$ such that $$u(x) \leq C \left ( d_g(x,
\overline{K_{\frac{\eta}{2}} \bigcup \{ x_1, ..., x_N \}}) \right )^{-\frac{2}{p-1}}.        $$

\end{itemize} 

\end{prop}

\noindent
This can be proved as in \cite{YY3} in the proof of Proposition 4.2, so we omit the details.

Now we are going to prove that $u$ is uniformly bounded on
$M \setminus K_{\eta}$. Suppose it is not, then there are sequences 
$\{ u_i \}$ and $ \{ p_i \} $ such that $$ \Delta _g u_i + Ku_i^{p_i}=0 \hspace{.3in} \text{ and } \hspace{.2in} \max _{M \setminus K_{\eta}} u_i
\rightarrow \infty \,\, \text{ as } i \rightarrow \infty.$$  
Therefore $\displaystyle \max _{M \setminus K_{\eta}}\left (
\left (d_g(x, K_{\frac{\eta}{2}}) \right )^{\frac{2}{p_i-1}}u_i(x)\right)\rightarrow \infty$
as $i \rightarrow \infty $. Then for fixed $\epsilon >0$ and $R >> 0
$  we can apply Proposition \ref{prop:iso} to each $u_i$ and
find $x_{1,i},..., x_{N(i),i}$ such that 
\begin{equation}
\label{locmax}
\text{ each } x_{j,i} \,\,
(1 \leq j \leq N(i) ) \text{ is a local maximum point of } u_i;
\end{equation}

\begin{equation}
\label{disjoint}
\text{the balls } B_{\frac{R}{u_i(x_{j, i})^{\frac{p_i -1}{2}}}}(x_{j, i}) \text{ are disjoint};   
\end{equation}

\noindent
for coordinates $y$ centered at $ x_{j,i}$ such that
$\frac{y}{u_i(x_{j,i})^{\frac{p_i-1}{2}}}$ is the conformally flat coordinate system,
\begin{equation}
\label{sphe}
\Bigg \| u_i(x_{j,i})^{-1}u_i \left (\frac
{y}{u_i(x_{j,i})^{\frac{p_i-1}{2}}}\right ) - \left (1+
\frac{K(x_{j, i})}{n(n-2)}|y|^2 \right )^{-\frac{n-2}{2}} \Bigg  \|_{C^2(B_{2R}(0))} <\epsilon;
\end{equation}
and
\begin{equation}
u_i(x) \leq C \left ( d_g(x,
\overline{K_{\frac{\eta}{2}} \bigcup \{ x_{1,i}, ..., x_{N(i),i} \}}) \right )^{-\frac{2}{p_i-1}} \hspace{.2in} \text{for a constant } C=C(\epsilon,R).       \end{equation}
 
Let $\sigma _i= \min \{ d_g(x_{\alpha, i}, x_{\beta ,i}): \alpha \neq
\beta, 1 \leq \alpha, \beta \leq N(i) \}$. Without lost of generality we
can assume $\sigma _i= d_g(x_{1, i}, x_{2 ,i})$. There are two
possibilities which could happen.

\noindent
{\bf Case I}: \hspace{.05in}$\sigma _i \geq \varepsilon >0$.\\Then the points $x_{j,i}$ have isolated
limiting points $x_1, x_2, ...$, which are isolated blow-up points of $\{
u_i \}$ as defined above.  

\noindent
{\bf Case II}: \hspace{.05in} $\sigma _i \rightarrow 0$.\\  Then we rescale the coordinates to make
the minimal distance 1: let $y=\sigma _i ^{-1}z$ where $z$ is the conformally flat coordinate system centered at $x_{1,i}$. We also rescale the function by defining
$$v_i(y)=\sigma _i ^{\frac{2}{p_i-1}}u_i(\sigma _i y).$$
$v_i$ satisfies 
$$
\Delta _{g^{(i)}} v_i + K(\sigma _i y)v_i^{p_i}=0
$$
where the metric
$g^{(i)}(y)=g_{\alpha \beta}(\sigma _i y)dy^{\alpha}dy^{\beta}$.
As proved in Section 4 of \cite{YY3}, $0$ is an isolated blow-up point of $\{ v_i \}$.

In Sections \ref {section:case1} and \ref {section:case2} we are going to prove that neither Case I nor Case
II can happen.

\section{Ruling out Case I}
\label{section:case1}

If the blow-up points are all isolated, then same argument as that in Section 6 of \cite{YY3} shows that among the isolated blow-up points $\{ x_1, x_2, ...\}$, there must be one which
is not a simple blow-up point, without loss of generality we assume it
to be $x_1$. To simplify the notations we are going to rename it to
be $x_0$. Let $ x_i $ be the local maximum point of $u_i$ such that
$\displaystyle \lim _{i \rightarrow \infty}x_i=x_0$.  

\noindent
Let $z=(z_1,...,z_n)$ be the conformally flat coordinates centered at each $x_i$. Since $x_0$ is not a simple blow-up point, as
a function of $|z|$, $|z|^{\frac{2}{p_i-1}}\bar{u}_i(|z|)$ has a
second critical point at $|z|=r_i$ where $r_i \rightarrow 0$. Let
$y=\frac{z}{r_i}$ and define
$v_i(y)=r_i^{\frac{2}{p_i-1}}u_i(r_iy)$. Then $v_i(y)$ satisfies
\begin{equation}
\label{eq:rescaled-case1}
\Delta _{g^{(i)}} v_i(y)+K_i(y)v_i(y)^{p_i}=0 
\end{equation} 
where $g^{(i)}(y) =g_{\alpha \beta}(r_iy)dy^{\alpha}dy^{\beta}$ and
$K_i(y)=K(r_iy)$.\\ 

\noindent
By this definition 
$|y|=1$ is the second critical point of $|y|^{\frac{2}{p_i-1}}\bar{v}_i(|y|)$. 
As shown in Section 6 of \cite{YY3}, $0$ is a simple blow-up point of $\{v_i\}$.

\subsection{Estimates for $v_i$}
\label {subsection:estforvi}

The following estimates are essentially the same as Proposition 5.3 in \cite{YY3}, except for a slightly different choice of parameters, but for completeness we repeat the proof.

\begin{prop}
\label{prop:simpleestimates}
There exists a constant $C$ independent of
$i$ such that
\begin{itemize}
\item if \,\, $0
\leq |y| \leq 1$, \hspace{.1in} then $$v_i(y) \geq C v_i(0)\left (1+ \frac{K_i(0)}{n(n-2)}v_i(0)^{\frac{4}{n-2}}|y|^2 \right )^{-\frac{n-2}{2}}$$ 
\item  if \,\,  $0 \leq |y| \leq \frac{R}{v_i(0)^{\frac{p_i-1}{2}}}$, \hspace{.1in}
then $$v_i(y) \leq C v_i(0) \left (1+ \frac{K_i(0)}{n(n-2)}
v_i(0)^{p_i-1}|y|^2 \right )^{-\frac{n-2}{2}}$$ 
\item if \,\, $\frac{R}{v_i(0)^{\frac{p_i-1}{2}}} \leq |y| \leq 1$,
\hspace{.1in} then \hspace{.1in} $v_i(y) \leq Cv_i(0)^{t _i}|y|^{-l_i}$\\
 where $l_i$, $t_i$ are so chosen that 
$\frac{2n-5}{2} < \displaystyle \lim_{i \to \infty} l_i < n-2$,
 and $t _i =1-\frac{(p_i-1)l_i}{2}$.
\end{itemize} 
\end{prop}

\pf
By Proposition \ref{prop:iso}, when $0 \leq |y| \leq \frac{R}{v_i(0)^{\frac{p_i-1}{2}}}$,
\begin{eqnarray*}
& & (1+\epsilon) v_i(0)\left (1+ \frac{K_i(0)}{n(n-2)}
v_i(0)^{p_i-1}|y|^2 \right )^{-\frac{n-2}{2}} \\
&  \geq & v_i(y)\\
& \geq  & (1-\epsilon) v_i(0)\left (1+ \frac{K_i(0)}{n(n-2)}
v_i(0)^{p_i-1}|y|^2 \right )^{-\frac{n-2}{2}}  \\
& \geq & (1-\epsilon) v_i(0)\left (1+ \frac{K_i(0)}{n(n-2)}
v_i(0)^{\frac{4}{n-2}}|y|^2\right )^{-\frac{n-2}{2}}. 
\end{eqnarray*}

\noindent
So we only need to find the upper and lower bounds on $v_i(y)$ when $\frac{R}{v_i(0)^{\frac{p_i-1}{2}}} \leq |y| \leq 1$.

\noindent
{\it First the lower bound.}

\noindent
Let $G_i$ be the Green's function of $\Delta _{g^{(i)}}$ which is singular at
$0$ and $G_i =0$ on $\partial B_1$. Since $\{g^{(i)}\}$
converges uniformly to the Euclidean metric, 
there exist constants $C_1$ and $C_2$ independent of $i$ such that
 $$C_1 |y|^{2-n} \leq G_i(y) \leq C_2 |y|^{2-n}.$$ 

\noindent
When $|y|=Rv_i(0)^{-\frac{p_i-1}{2}}$,
\begin{eqnarray*}
v_i(y) & \geq & (1-\epsilon) \frac{v_i(0)}{\left (1+ \frac{K_i(0)}{n(n-2)}
v_i(0)^{p_i-1}|y|^2 \right )^{\frac{n-2}{2}}}\\
& = & (1-\epsilon)\frac{v_i(0)}{\left (1+ \frac{K_i(0)}{n(n-2)}
R^2 \right )^{\frac{n-2}{2}}}\\
& = & (1-\epsilon)\left( R^{-2}+\frac{K_i(0)}{n(n-2)}\right)^{-\frac{n-2}{2}}R^{2-n}v_i(0)\\
& \geq & CR^{2-n}v_i(0)\\
& \geq & CR^{2-n}v_i(0)^{\frac{(n-2)(p_i-1)}{2}-1} \hspace{.2in} \text{ since }\hspace{.1in} \frac{(n-2)(p_i-1)}{2}-1 \leq 1\\ 
& = & Cv_i(0)^{-1}|y|^{2-n}\\
& \geq & Cv_i(0)^{-1} G_i(y)
\end{eqnarray*}

\noindent
With this constant $C$, when $|y|=1$, $Cv_i(0)^{-1} G_i(y)=0<v_i(y)$.

\noindent
We know that $$\Delta _{g^{(i)}}\left (v_i(y)- Cv_i(0)^{-1}G_i(y)\right )=\Delta_{g^{(i)}}v_i(y)=-K_i(y)v_i(y)^{p_i}<0$$ on $B_1 \setminus B_{Rv_i(0)^{-\frac{p_i-1}{2}}}$.  Therefore, by the maximal principle, when $\frac{R}{v_i(0)^{\frac{p_i-1}{2}}} \leq |y| \leq 1$,
\begin{eqnarray*}
 v_i(y) & > & Cv_i(0)^{-1}G_i(y)\\
 & \geq & Cv_i(0)^{-1}|y|^{2-n}.
\end{eqnarray*}

\noindent
Now we need to compare $|y|^{2-n}v_i(0)^{-1}$ with $v_i(0)\cdot\left (1+ \frac{K_i(0)}{n(n-2)}v_i(0)^{\frac{4}{n-2}}|y|^2 \right )^{-\frac{n-2}{2}}$ in order to get the desired lower bound.

\begin{eqnarray*}
 & & v_i(0)^2|y|^{n-2} \left (1+ \frac{K_i(0)}{n(n-2)}
v_i(0)^{\frac{4}{n-2}}|y|^2 \right )^{-\frac{n-2}{2}} \\
& \leq &  v_i(0)^2 \left (\frac{K_i(0)}{n(n-2)}
v_i(0)^{\frac{4}{n-2}} \right )^{-\frac{n-2}{2}}\\
& \leq & C  
\end{eqnarray*} 
for a constant $C$ independent of $i$. Therefore $$v_i(0)^{-1}|y|^{2-n}
\geq C v_i(0) \left (1+ \frac{K_i(0)}{n(n-2)}
v_i(0)^{\frac{4}{n-2}}|y|^2 \right )^{-\frac{n-2}{2}}, $$
and consequently
$$v_i(y) \geq C v_i(0) \left (1+ \frac{K_i(0)}{n(n-2)}
v_i(0)^{\frac{4}{n-2}}|y|^2 \right )^{-\frac{n-2}{2}}$$ 
when $\frac{R}{v_i(0)^{\frac{p_i-1}{2}}} \leq |y| \leq 1$.

\noindent
{\it Next the upper bound.} 

\noindent
We are going to apply the same strategy of
constructing a comparison function and using the maximal principle.

\noindent
Define $\mathcal{L}_i\varphi := \Delta _{g^{(i)}} \varphi + K_i
v_i^{p_i-1}\varphi$. By this definition $\mathcal{L}_iv_i=0$. Let $M_i=\displaystyle \max _{\partial B_1} v_i$ and $C_i=(1+\epsilon)
\left (\frac{K_i(0)}{n(n-2)} \right )^{-\frac{n-2}{2}}$. Note that $C_i$ is bounded above
and below by constants independent of $i$. Consider the function
$$M_i|y|^{-n+2+l_i}+C_iv_i(0)^{t _i}|y|^{-l_i}.$$

\noindent
When $|y|=\frac{R}{v_i(0)^{\frac{p_i-1}{2}}}$,
 
\begin{eqnarray*}
v_i(y) & \leq & (1+\epsilon) \frac{v_i(0)}{\left (1+ \frac{K_i(0)}{n(n-2)}
v_i(0)^{p_i-1}|y|^2 \right )^{\frac{n-2}{2}}} \\
&=& (1+\epsilon) \frac{v_i(0)}{\left (1+ \frac{K_i(0)}{n(n-2)}R^2
\right )^{\frac{n-2}{2}}} \\
& \leq & C_iv_i(0)R^{-(n-2)}\\
& \leq & C_iv_i(0)R^{-l_i}  \\
& = & C_i v_i(0)^{t _i}|y|^{-l _i}. 
\end{eqnarray*}

\noindent
When $|y|=1$, by the definition of $M_i$, $v_i(y) \leq M_i=M_i|y|^{-n+2+l_i}.$ \\
Thus on $\{|y|=1 \} \cup \{|y|=Rv_i(0)^{-\frac{p_i-1}{2}}\}$,$$v_i(y) \leq M_i|y|^{-n+2+l_i}+C_iv_i(0)^{t _i}|y|^{-l_i}.$$

\noindent
In the Euclidean coordinates, $\Delta |y|^{-l_i}=-l_i(n-2-l_i)|y|^{-l_i-2}
$ and $\Delta |y|^{-n+2+l_i}=-l_i(n-2-l_i)|y|^{-n+l_i}$.  When $i$ is sufficiently large, $g^{(i)}$ is close to the Euclidean metric. Therefore

\begin{equation}
\label{eq:|y|^-l_i}
\Delta _{g^{(i)}} |y|^{-l_i} \leq -\frac{1}{2}l_i(n-2-l_i)|y|^{-l_i-2}
\end{equation}
and 
\begin{equation}
\label{eq:|y|^-n+2+l_i}
\Delta _{g^{(i)}} |y|^{-n+2+l_i} \leq
-\frac{1}{2}l_i(n-2-l_i)|y|^{-n+l_i}. 
\end{equation}
Thus
\begin{eqnarray*}
& & \mathcal{L}_i(C_iv_i(0)^{t _i}|y|^{-l_i}) \\
&=& C_i v_i(0)^{t_i}\Delta _{g^{(i)}}|y|^{-l_i}+C_iv_i(0)^{t _i}K_iv_i(y)^{p_i-1}|y|^{-l_i} \\
& \leq & -Cl_i(n-2-l_i)v_i(0)^{t _i}|y|^{-l_i-2}+C'v_i(0)^{t _i}v_i(y)^{p_i-1}|y|^{-l_i}
\end{eqnarray*}
for some constants $C, C'$ independent of $i$.

\noindent
Lemma \ref{lemma:simpleharnack} and the upper bound on $v_i(y)$ when $|y|\leq Rv_i(0)^{-\frac{p_i-1}{2}}$ imply
that
\begin{eqnarray*}
\bar{v}_i \left (Rv_i(0)^{-\frac{p_i-1}{2}} \right ) & \leq & 
C\frac{(1+\epsilon)v_i(0)}{ \left [1+ \frac{K_i(0)}{n(n-2)}
v_i(0)^{p_i-1} \left ( Rv_i(0)^{-\frac{p_i-1}{2}} \right )^2 \right ]^{\frac{n-2}{2}}}\\
& \leq & Cv_i(0)R^{2-n}.
\end{eqnarray*}

\noindent
Then since $0$ is a simple blow-up point and $r^{\frac{2}{p_i-1}}\bar{v}
_i(r)$ is decreasing from $Rv_i(0)^{-\frac{p_i-1}{2}}$ to $1$,

\begin{eqnarray*}
|y|^{\frac{2}{p_i-1}} \bar{v}_i(|y|) & \leq &
\left (Rv_i(0)^{-\frac{p_i-1}{2}} \right )^{\frac{2}{p_i-1}} \cdot \bar{v}_i
\left (Rv_i(0)^{-\frac{p_i-1}{2}} \right )\\
& \leq & CR^{\frac{2}{p_i-1}+2-n}.
\end{eqnarray*}
Thus again by Lemma \ref{lemma:simpleharnack}  
\begin{equation}
\label{eq:v_i^p_i-1|y|^-l_i}
v_i(y)^{p_i-1} \,\, \leq C\bar{v}_i(|y|)^{p_i-1} \,\, \leq \,\, C|y|^{-2}R^{2-(n-2)(p_i-1)},
\end{equation}
and hence 
$$v_i(y)^{p_i-1} |y|^{-l_i}\leq C|y|^{-2-l_i}R^{2-(n-2)(p_i-1)}.$$
Therefore
\begin{eqnarray*}
& & \mathcal{L}_i \left (C_iv_i(0)^{t _i}|y|^{-l_i} \right )\\
 &   \leq & \left(-Cl_i(n-2-l_i)+C' R^{2-(n-2)(p_i-1)} \right )  v_i(0)^{t _i}|y|^{-l_i-2} 
\end{eqnarray*}
By our choice of $l_i$, \,\, $l_i(n-2-l_i)$ is always
bounded below by some positive constant independent of $i$. When $i$
is sufficiently large, $2-(n-2)(p_i-1) < 0$, so we can choose $R$ big enough
such that $-Cl_i(n-2-l_i)+C' R^{2-(n-2)(p_i-1)}<0 $, which implies $\mathcal{L}_i(C_iv_i(0)^{t _i}|y|^{-l_i}) <0$.

\noindent
Similarly,

\begin{eqnarray*}
\mathcal{L}_i \left (M_i|y|^{-n+2+l_i} \right ) & = & M_i \Delta _{g^{(i)}} |y|^{-n +2 + l_i} +M_iK_iv_i^{p_i-1}|y|^{-n+2+l_i}\\
& \leq & 
-\frac{1}{2}l_i(n-2-l_i)M_i|y|^{-n+l_i}
+K_i M_i R^{2-(n-2)(p_i-1)}|y|^{-n+l_i} 
\end{eqnarray*}
by equations (\ref{eq:|y|^-n+2+l_i}) and (\ref{eq:v_i^p_i-1|y|^-l_i}). We can choose $R$ large enough such that
$-\frac{1}{2}l_i(n-2-l_i) +K_iR^{2-(n-2)(p_i-1)} <0$ and hence
$$\mathcal{L}_i(M_i|y|^{-n+2+l_i})<0. $$

\noindent
Therefore when $Rv_i(0)^{-\frac{p_i-1}{2}} \leq |y| \leq 1$,
$$\mathcal{L}_i \left (M_i|y|^{-n+2+l_i}+C_iv_i(0)^{t
_i}|y|^{-l_i} \right ) <0. $$  Then by the maximal principle  
$$
v_i(y)  \leq  M_i|y|^{-n+2+l_i}+C_iv_i(0)^{t_i}|y|^{-l_i}.
$$

\noindent
By Lemma \ref{lemma:simpleharnack} and because $0$ is a simple blow-up point, for
$\frac{R}{v_i(0)^{\frac{p_i-1}{2}}} \leq   \theta \leq 1$, 
\begin{eqnarray*}
M_i & \leq & C \theta^{\frac{2}{p_i-1}}\bar{v}_i(\theta) \\
& \leq & C\theta ^{\frac{2}{p_i-1}} \left (M_i\theta^{-n+2+l_i}+C_iv_i(0)^{t
_i} \theta ^{-l_i} \right )\\
 & = & C \theta ^{\frac{2}{p_i-1}-n+2+l_i}M_i + C \theta
^{\frac{2}{p_i-1}} \cdot C_iv_i(0)^{t_i} \theta ^{-l_i} 
\end{eqnarray*}
for some constant $C$ independent of $i$.

\noindent
Note that $$ \displaystyle \lim_{i \to \infty}\left(\frac{2}{p_i-1}-n+2+l_i \right)\,\, = \,\, -\frac{n-2}{2}+  \displaystyle \lim_{i \rightarrow \infty} l_i \,\, > \,\,  -\frac{n-2}{2}+ \frac{2n-5}{2} \,\, > \,\, 0$$ because $n \geq 5$.

\noindent
Since $\frac{R}{v_i(0)^{\frac{p_i-1}{2}}} \to 0$, we can choose $\theta$ small enough (fixed and independent of $i$) to
absorb the first term on the right hand side of the above inequality
into the left hand side to get $$M_i \,\, \leq \,\, 2C \theta
^{\frac{2}{p_i-1}}\cdot C_i v_i(0)^{t_i} \theta ^{-l_i}\,\, 
\leq \,\, C v_i(0)^{t _i}.$$ 

\noindent
Therefore
\begin{eqnarray*}
v_i(y) & \leq &  M_i|y|^{-n+2+l_i}+C_iv_i(0)^{t_i}|y|^{-l_i}\\
& \leq & M_i|y|^{-l_i}+C_iv_i(0)^{t_i}|y|^{-l_i} 
\\
& \leq &  C v_i(0)^{t _i}|y|^{-l_i}
\end{eqnarray*}

\stop

\subsection{A Preliminary Estimate for $\delta _i :=\frac {n+2}{n-2} - p_i$}
\label{subsection:prelim-delta}

First we prove a technical lemma.

\begin{lemma}
\label{lemma:intv^p+1y^k}
When $\sigma <1$ and $0 \leq \kappa \leq n-2$, $$\int _{|y|\leq \sigma} |y|^{\kappa}v_i(y)^{p_i+1}dy \leq C v_i(0)^{-\frac{2\kappa}{n-2}+\frac{n-2+\kappa}{2}\delta_i},$$ where $C$ is independent of $i$.
\end{lemma}

\pf
\noindent
By Proposition \ref{prop:simpleestimates}
\begin{eqnarray*}
 \int _{|y| \leq \frac{R}{v_i(0)^{\frac{p_i-1}{2}}} } |y|^{\kappa}
v_i(y)^{p_i+1} dy & \leq & C
v_i(0)^{p_i+1}  \int _{|y| \leq \frac{R}{v_i(0)^{\frac{p_i-1}{2}}}
} |y|^{\kappa} dy   \\
& \leq & C v_i(0)^{p_i+1-\frac{(n+\kappa)(p_i-1)}{2}}\\
& = & C v_i(0) ^{-\frac{2\kappa}{n-2} + \frac{n-2+\kappa}{2}\delta _i}.
\end{eqnarray*}

\noindent
Since $n \geq 5$, by our choice of $\l_i$ 
\begin{eqnarray*}
\lim _{i \rightarrow \infty} \big ( n+\kappa-l_i(p_i+1) \big ) &  = &  n+\kappa-\frac{2n}{n-2}\lim_{i \to \infty} l_i \\
& < & n+ \kappa-\frac{2n}{n-2}\cdot \frac{2n-5}{2}\\
& \leq & n+(n-2)-\frac{n(2n-5)}{n-2}\\
& < & 0.
\end{eqnarray*}
Therefore
\begin{eqnarray*}
\int _{\frac{R}{v_i(0)^{\frac{p_i-1}{2}}} \leq |y| \leq \sigma } |y|^{\kappa}v_i(y)^{p_i+1} dy & \leq & C \int
_{\frac{R}{v_i(0)^{\frac{p_i-1}{2}}} \leq |y| \leq \sigma } |y|^{\kappa}
\left ( v_i(0)^{t _i}|y|^{-l_i} \right )^{p_i+1} dy \\
            & \leq & C v_i(0)^{ t _i(p_i+1) -\frac{p_i-1}{2}\left
(n-l_i(p_i+1)+\kappa\right )} \\
& = & C v_i(0)^{p_i+1-\frac{(n+\kappa)(p_i-1)}{2}} \hspace{.2in}(\text{by the definition of }  t _i)\\
 & = & C v_i(0) ^{-\frac{2\kappa}{n-2} +\frac{n-2+\kappa}{2}\delta_i}.
\end{eqnarray*}
Thus
\begin{equation*}
\int _{|y| \leq \sigma } |y|^{\kappa} v_i(y)^{p_i+1} dy \leq C v_i(0)
^{-\frac{2\kappa}{n-2} + \frac{n-2+\kappa}{2}\delta _i}.
\end{equation*}

\stop

The next proposition is a preliminary estimate for $\delta _i:=  \frac {n+2}{n-2} - p_i $, we will also derive a refined estimate in a later part of this paper.

\begin{prop}
\label{prop:v^delta}
$\displaystyle \lim _{i \rightarrow \infty} v_i(0)^{\delta _i} =1$.
\end{prop}

\pf 
Since the original metric is locally conformally flat, locally it can be written as $\lambda (z) ^{\frac{4}{n-2}} dz^2. $ Let $\lambda _i(y)=\lambda(r_iy)$, then $g^{(i)}(y)=\lambda _i(y)^{\frac{4}{n-2}}dy^2 $.  Let $\sigma < 1 $,
the Pohozaev identity in \cite{S5} says that for a conformal Killing field $X$ on $B_{\sigma}$,  
\begin{equation}
\label{eq:pohozaev-X(R)} 
\frac{n-2}{2n} \int _{B_{\sigma}} X(R_i) dv_{g_i} = \int _{\partial
B_{\sigma}} T_i(X, \nu _i) d \Sigma _i
\end{equation}

\noindent
where the notations are
\begin{eqnarray*}
g_i &  = &  v_i^{\frac{4}{n-2}}g^{(i)} \,\,  = \,\, (\lambda_i v_i)^{\frac{4}{n-2}}dy^2, \\
R_i &  = & R(g_i) \,\, = \,\, c(n)^{-1} K_iv_i^{- \delta _i},\\
 dv_{g_i} & = &  (\lambda _i v_i)^{\frac{2n}{n-2}}dy  ,\\
\nu _i & = & (\lambda _i v_i)^{-\frac{2}{n-2}}\sigma^{-1}\sum _{j} y^j
\frac{\partial}{\partial y^j} \\
& &  \text { is the unit outer
normal vector on } \partial B_{\sigma} \text{ with respect to } g_i,\\
d\Sigma _i & = & (\lambda _i v_i)^{\frac{2(n-1)}{n-2}} d \Sigma _{\sigma}
\\
& & 
\text{ where } d \Sigma _{\sigma} \text { is the surface element of the
standard } S^{n-1}(\sigma),\\
T_i & = & \Ric(g_i)-n^{-1}R(g_i)g_i \hspace{.2in} \text{ is the traceless Ricci tensor with respect to } g_i.
\end{eqnarray*}

\noindent
$T_i$ can also be expressed as  (see \cite{S7}) $$(n-2) (\lambda _i v_i)^{\frac{2}{n-2}} \left(\Hess \left ( (\lambda _i v_i)^{-\frac{2}{n-2}} \right )-\frac{1}{n}\Delta \left ((\lambda _i
v_i)^{-\frac{2}{n-2}} \right ) dy^2 \right )$$
\noindent
where $\Hess$ and $\Delta$ are taken with respect to the Euclidean metric $dy^2$.

We choose $X=\displaystyle \sum _{j=1}^{n} y^j \frac{\partial}{\partial y^j}$.

\noindent 
The left hand side is
\begin{eqnarray*}
 & & \frac{n-2}{2n}\int _{B_{\sigma}} X(R_i) dv_{g_i} \\
& = &\frac{2(n-1)}{n}\int _{B_{\sigma}} X(K_iv_i^{-\delta _i}) (\lambda _i v_i)^{\frac{2n}{n-2}} dy \\
& = & \frac{2(n-1)}{n} \int _{B_{\sigma}} X(K_i)v_i^{p_i+1} \lambda _i ^{\frac{2n}{n-2}} dy - \frac{2(n-1)}{n}\delta _i \int  _{B_{\sigma}} K_i v_i^{p_i}X(v_i)
\lambda _i^{\frac{2n}{n-2}} dy  
\end{eqnarray*}
By the divergence theorem it is equal to 
\begin{eqnarray*}
& = & \frac{2(n-1)}{n} \int  _{B_{\sigma}}  |y| \frac {\partial K_i}{\partial r} v_i^{p_i+1} \lambda _i ^{\frac{2n}{n-2}}dy +
\frac{2(n-1)}{n}\frac{\delta _i}{p_i+1} \bigg (\int _{B_{\sigma}} r \frac{ \partial
K_i}{\partial r} \lambda_i ^{\frac{2n}{n-2}} v_i^{p_i+1} dy  \\
& & +  \int_{B_{\sigma}} K_i v_i^{p_i+1} r \frac{ \partial \lambda _i^{\frac{2n}{n-2}}}{\partial r} dy  + \int
_{ B_{\sigma}} K_i v_i^{p_i+1}\lambda _i^{\frac{2n}{n-2}} \diver X \,\, dy \bigg )\\
& &  -  \frac{2(n-1)}{n}\frac{\delta _i}{p_i+1}   \int
_{\partial B_{\sigma}} K_iv_i^{p_i+1}\lambda _i^{\frac{2n}{n-2}}  X \cdot
\left (\frac{\sum
y^j \frac{\partial}{\partial y^j}}{\sigma} \right )d\Sigma _{\sigma},
\end{eqnarray*}
which can be further written as

\begin{eqnarray}
\label{eq:X(R)in-v^delta}
& = & \frac{2(n-1)}{n}\left(1+ \frac{\delta _i}{p_i+1}\right) \int  _{B_{\sigma}}  |y| \frac {\partial K_i}{\partial r} v_i^{p_i+1}\lambda _i^{\frac{2n}{n-2}}dy  \nonumber \\
& &  + \frac{2(n-1)}{n}\frac{\delta _i}{p_i+1} \int
_{B_{\sigma}}|y| K_iv_i^{p_i+1}  \frac{ \partial \lambda _i^{\frac{2n}{n-2}}}{\partial r} dy  \\
& &  +\frac{2(n-1)}{n}\frac{\delta _i}{p_i+1}n  \int
_{ B_{\sigma}} K_iv_i^{p_i+1}\lambda _i^{\frac{2n}{n-2}} dy -\frac{2(n-1)}{n} \frac{\delta _i}{p_i+1} \int _{\partial B_{\sigma}} \sigma K_iv_i^{p_i+1}\lambda_i ^{\frac{2n}{n-2}}d\Sigma _{\sigma}. \nonumber
\end{eqnarray}

The right hand side of (\ref{eq:pohozaev-X(R)}) is
\begin{eqnarray}
\label{eq:v^delta-Tterm}
& & \int _{\partial
B_{\sigma}} T_i(X, \nu _i) d \Sigma _i \nonumber \\
 & = & \int _{\partial
B_{\sigma}} (n-2) (\lambda _i v_i)^{\frac{2}{n-2}}
\bigg [ \Hess \left ((\lambda _i v_i)^{-\frac{2}{n-2}} \right ) 
  \left (r \frac {\partial}{\partial r}, (\lambda _iv_i)^{-\frac{2}{n-2}}\sigma^{-1} r \frac {\partial}{\partial r}
  \right )  \nonumber \\  
& &    -\frac{1}{n}\Delta 
  \left ((\lambda _i v_i)^{-\frac{2}{n-2}} \right) 
  \left < r \frac {\partial}{\partial r},  (\lambda _i
  v_i)^{-\frac{2}{n-2}}\sigma^{-1} r \frac {\partial}{\partial r} \right > 
\bigg ] 
(\lambda _i v_i)^{\frac{2(n-1)}{n-2}}  d \Sigma _{\sigma} \nonumber \\
& & (\text{where } <\cdot, \cdot> \text { is the Euclidean metric})
\nonumber \\
& = & (n-2) \int _{\partial B_{\sigma}} \bigg [\sigma ^{-1}\Hess 
\left ( (\lambda _i v_i)^{-\frac{2}{n-2}} \right ) 
\left  (r \frac {\partial}{\partial r}, r \frac {\partial}{\partial r}
\right ) \nonumber \\
& & -\frac{\sigma}{n}\Delta \left ((\lambda _i
v_i)^{-\frac{2}{n-2}} \right ) \bigg ] (\lambda _i
v_i)^{\frac{2(n-1)}{n-2}} d \Sigma _{\sigma} \nonumber \\
& = & (n-2)
\int _{\partial B_{\sigma}} 
 \sigma ^{-1}  \bigg [ -\frac {2 }{n-2} (\lambda _i v_i) \sum _{j,k}y^j y^k 
\frac{\partial}{\partial y^k} \frac{\partial}{\partial y^j} (\lambda _i
v_i)  \nonumber \\
& & +  \frac{2n}{(n-2)^2} \sum_{j,k} y^j y^k
\frac{\partial (\lambda _i v_i)}{\partial y^k} \frac{\partial( \lambda _i
v_i)}{\partial y^j}  \bigg ] \\ 
& & - \sigma \cdot \bigg [ -\frac{2 }{n(n-2)} (\lambda _i v_i) \sum _{j} \frac{\partial ^2 (\lambda _i v_i)}{(\partial
y^j)^2}  + \frac{2}{(n-2)^2} \sum _{j} \left(\frac{\partial (\lambda _i
v_i)}{\partial y^j} \right )^2 \bigg ]
d \Sigma _{\sigma} \nonumber
\end{eqnarray}

Next we are going to study the decay rate of each term in (\ref {eq:X(R)in-v^delta}) and (\ref {eq:v^delta-Tterm}). 

\noindent
On $\partial B _{\sigma} $, by Proposition \ref{prop:simpleestimates}, $v_i \leq C
v_i(0)^{t _i}$, then by the elliptic regularity theory \cite{GT} $\| v_i \| _{C^2(\partial B _{\sigma}) } \leq C
v_i(0)^{t _i}$. Thus we know (\ref{eq:v^delta-Tterm}) decays in the rate of $v_i(0)^{2t _i}$. 

\noindent
The fourth term in (\ref{eq:X(R)in-v^delta}) decays in
the order of $\delta _iv_i(0)^{t _i(p_i+1)}$ by Proposition \ref{prop:simpleestimates}.  By Lemma \ref {lemma:intv^p+1y^k} we know that the second term in (\ref{eq:X(R)in-v^delta}) is bounded above by $$C\delta_i \int _{B_{\sigma}}|y|v_i^{p_i+1}dy \leq C\delta_i v_i(0)^{-\frac{2}{n-2}+\frac{n-1}{2}\delta_i}.$$

\noindent
Therefore the sum of the first and the third terms in (\ref{eq:X(R)in-v^delta}), which is $$\frac{n}{2(n-1)}\left(1+ \frac{\delta _i}{p_i+1}\right) \int  _{B_{\sigma}}  |y| \frac {\partial K_i}{\partial r} v_i^{p_i+1}\lambda _i^{\frac{2n}{n-2}}dy  + \frac{n}{2(n-1)}\frac{\delta _i}{p_i+1}n  \int
_{ B_{\sigma}} K_iv_i^{p_i+1}\lambda _i^{\frac{2n}{n-2}} dy$$ is bounded above by $Cv_i(0)^{2t_i}+C\delta _i v_i(0)^{t_i(p_i+1)} + C\delta_i v_i(0)^{-\frac{2}{n-2}+\frac{n-1}{2}\delta_i} $.

\noindent
By our choice of $l_i$ and $t_i$, as $i \to \infty$,  $$  t _i = 1 -\frac{(p_i-1)l_i}{2}\,\, \rightarrow \,\, 1-\frac{2}{n-2} \lim _{i \to \infty} l_i \,\, < \,\,  1-\frac{2}{n-2} \cdot \frac{2n-5}{2} \,\, < \,\, 0.$$
Thus $Cv_i(0)^{2t_i}+C\delta _i v_i(0)^{t_i(p_i+1)} \leq Cv_i(0)^{2t_i}+Cv_i(0)^{t_i(p_i+1)} \leq Cv_i(0)^{2t_i} $.

\noindent
On the other hand 
$$\frac{\delta _i}{p_i+1}n  \int
_{ B_{\sigma}} K_iv_i^{p_i+1}\lambda _i^{\frac{2n}{n-2}} dy \geq C \delta
_i   \int _{ B_{\sigma}} v_i^{p_i+1} dy .$$
When $|y| \leq \frac{R}{v_i(0)^{\frac{p_i-1}{2}}}$, by Proposition \ref{prop:simpleestimates} 
\begin{eqnarray*}
 v_i(y)& \geq & (1-\epsilon) \frac{v_i(0)}{\left (1+ \frac{K_i(0)}{n(n-2)}
v_i(0)^{p_i-1}|y|^2 \right )^{\frac{n-2}{2}}}\\
& \geq & (1-\epsilon)
\frac{v_i(0)}{\left (1+ \frac{K_i(0)}{n(n-2)}R^2 \right )^{\frac{n-2}{2}}}\\
& \geq &  Cv_i(0),
\end{eqnarray*}
so
\begin{eqnarray}
\label{eq:lower-int-v^p+1}
 \int _{ B_{\sigma}} v_i^{p_i+1} dy & > & \int _{|y| \leq
\frac{R}{v_i(0)^{\frac{p_i-1}{2}} } } v_i^{p_i+1} dy \nonumber \\
& \geq & C v_i(0)^{p_i+1 - \frac{n}{2}(p_i-1)} \nonumber \\
& = & C v_i(0)^{ \frac{n-2}{2} \delta _i} \nonumber \\
& \geq & C.
\end{eqnarray}
This implies that the third term in (\ref{eq:X(R)in-v^delta}) is bounded below by $C \delta _i$.

\noindent
Then by comparing the decay rates of the terms in (\ref {eq:X(R)in-v^delta}) and (\ref {eq:v^delta-Tterm}), 
$$\delta _i \leq C \left ( v_i(0)^{2t _i} + \delta_i v_i(0)^{-\frac{2}{n-2}+\frac{n-1}{2}\delta_i}  + \bigg | \int  _{B_{\sigma}}   \frac {\partial K_i}{\partial r}|y| v_i^{p_i+1}\lambda_i^{\frac{2n}{n-2}}dy \bigg | \right ).$$  Since $v_i(0)^{-\frac{2}{n-2}+\frac{n-1}{2}\delta_i} \to 0$, the second term on the right hand side can be absorbed into the left hand side.
Thus we conclude that 
\begin{equation}
\label{eq:deltacoarse}
\delta _i \leq C \left ( v_i(0)^{2t _i} + \bigg |\int  _{B_{\sigma}}   \frac {\partial K_i}{\partial r}|y| v_i^{p_i+1}\lambda_i^{\frac{2n}{n-2}}dy  \bigg|\right ) .
\end{equation} 

\noindent
By Lemma \ref {lemma:intv^p+1y^k}, $ \big | \int  _{B_{\sigma}}   \frac {\partial K_i}{\partial r}|y| v_i^{p_i+1}\lambda_i^{\frac{2n}{n-2}}dy \big | \leq C v_i(0)^{-\frac{2}{n-2}+\frac{n-1}{2}\delta_i}$, thus 

$$
 \delta _i \leq C\left (v_i(0)^{-\frac{2}{n-2}+\frac{n-1}{2}\delta _i} +
v_i(0)^{2t _i} \right ).
$$

\noindent
This implies that 
$$
 \delta _i \ln v_i(0)\leq
C\left (v_i(0)^{-\frac{2}{n-2}+\frac{n-1}{2}\delta _i} +
v_i(0)^{2t _i} \right )\ln v_i(0)  \rightarrow  0$$ as $i \to \infty$.  Therefore $\displaystyle \lim_{i \to \infty} v_i(0)^{\delta_i} =1$.  Consequently, we have

\begin{equation}
\label{eq:deltaiprelim}
 \delta _i \leq C \left ( v_i(0)^{-\frac{2}{n-2}}+v_i(0)^{2t _i} \right ).
\end{equation}

\stop

\subsection{A Preliminary Estimate for $|\nabla K_i|$}
\label{subsection:prelim-gradK}

We will again study the Pohozaev identity (\ref{eq:pohozaev-X(R)}), but with a different choice of the conformal Killing field $X=\frac
{\partial}{\partial y^1}$. 

\noindent
Direct calculation, as that in the proof of Proposition \ref{prop:v^delta}, shows that the right hand side of the identity is equal to 
\begin{eqnarray*}
& & (n-2)
\int _{\partial B_{\sigma}} 
 \sum _{j} \frac{y^j}{\sigma}    \left( -\frac {2 }{n-2} (\lambda _i v_i) 
\frac{\partial ^2 (\lambda _i v_i)  }{\partial y^1 \partial y^j} 
 +  \frac{2n}{(n-2)^2}
\frac{\partial (\lambda _i v_i)}{\partial y^1} \frac{\partial( \lambda _i
v_i)}{\partial y^j}  \right) \\ 
& & - \frac{y^1}{\sigma} \sum _{j}\left( -\frac{2 }{n(n-2)} (\lambda _i v_i)  \frac{\partial ^2 (\lambda _i v_i) }{(\partial y^j)^2} + \frac{2}{(n-2)^2} \left(\frac{\partial (\lambda _i
v_i)}{\partial y^j} \right )^2 \right )
d \Sigma _{\sigma}, 
\end{eqnarray*}
and decays in the rate of $ v_i(0)^{2 t _i}$.

\noindent
The left hand side of this identity is 
\begin{eqnarray}
\label{eq:derivativeofK-lhs}
& &  \frac{n-2}{2n} \int _{B_{\sigma}} \frac{\partial}{\partial
y^1}(R_i) dv_{g_i} \nonumber \\
& = & \frac{n-2}{2n} c(n)^{-1} \int _{B_{\sigma}} \frac{\partial}{\partial
y^1}(K_iv_i^{-\delta _i} )(\lambda_i v_i)^{\frac{2n}{n-2}}dy \nonumber \\
& = & \frac{n-2}{2n} c(n)^{-1} \int _{B_{\sigma}} \left (1+
\frac{\delta_i}{p_i+1} \right )\lambda_i^
{\frac{2n}{n-2}}v_i^{p_i+1}\frac{\partial K_i}{\partial
y^1} dy \nonumber \\
& & +\frac{n-2}{2n}c(n)^{-1} \int _{B_{\sigma}} \frac{\delta_i}{p_i+1}
K_i v_i^{p_i+1} \frac {\partial \lambda _i^
{\frac{2n}{n-2}}}{\partial y^1} dy  \nonumber\\
& & - \frac{n-2}{2n} c(n)^{-1}
\frac{\delta_i}{p_i+1}\int _{\partial B_{\sigma}}
\lambda_i^{\frac{2n}{n-2}}K_iv_i^{p_i+1} \frac{y^1}{\sigma}
d \Sigma _{\sigma}.
\end{eqnarray}

\noindent
By Proposition \ref{prop:simpleestimates}, the last term in (\ref{eq:derivativeofK-lhs}) is
bounded above by $$C\delta _i \cdot v_i(0)^{t_i(p_i+1)} \leq C\delta _iv_i(0)^{2t _i}
$$ since $t_i < 0$ and $v_i(0) \rightarrow \infty$.

\noindent
Note that 
$\lambda _i(y)=\lambda (r_iy)$, the second term in (\ref{eq:derivativeofK-lhs}) is bounded above by 
$$C\delta _i r_i \int_{|y|\leq \sigma} v_i(y)^{p_i+1}dy, $$ which is further bounded by $C\delta_i r_iv_i(0)^{\frac{n-2}{2}\delta_i}\leq C\delta _i r_i$ by Lemma \ref {lemma:intv^p+1y^k} and Proposition \ref{prop:v^delta}.

\noindent
Therefore the first term in (\ref{eq:derivativeofK-lhs}) which is $$\frac{n-2}{2n} c(n)^{-1} \int _{B_{\sigma}} \left(1+
\frac{\delta_i}{p_i+1}\right)\lambda_i^
{\frac{2n}{n-2}}v_i^{p_i+1}\frac{\partial K_i}{\partial
y^1} dy$$ is bounded above by 
$
C(v_i(0)^{2 t _i}+\delta _i v_i(0)^{2t_i}+ \delta _i r_i) \leq C \left ( \delta _i r_i + v_i(0)^{2 t _i} \right )$.

\noindent
This shows that 
\begin{equation}
\label{eq:gradKicoarse}
\bigg| \int _{B_{\sigma}} \lambda _i^{\frac{2n}{n-2}}v_i^{p_i+1}\frac{\partial K_i}{\partial
y^1} dy \bigg| \leq C \left ( \delta _i r_i + v_i(0)^{2 t _i} \right ).
\end{equation}

\noindent
By the Taylor expansion
$$\frac{\partial K_i}{\partial y^1}(y)= \frac{\partial K_i}{\partial
y^1}(0) + \nabla \left(\frac{\partial K_i}{\partial
y^1} \right )(\varsigma) \cdot y
\hspace{.3in} \text{ for some } |\varsigma| \leq |y|.$$

\noindent
Note that $K _i(y)=K(r_iy)$.  By Lemma \ref{lemma:intv^p+1y^k} and Proposition \ref{prop:v^delta},

\begin{eqnarray*}
 \int _{B_{\sigma}} \lambda _i^{\frac{2n}{n-2}}v_i^{p_i+1} \Bigg | \nabla \left (\frac{\partial K_i}{\partial y^1} \right )(\varsigma) \cdot y  \Bigg |dy & \leq & C r_i \int _{B_{\sigma}}
v_i^{p_i+1}|y|dy \\
& \leq & C r_i v_i(0)
^{-\frac{2}{n-2} +\frac{n-1}{2}\delta_i}  \\
& \leq & Cr_iv_i(0)^{-\frac{2}{n-2}}.  
\end{eqnarray*}

\noindent
Thus we know 
\begin{eqnarray*}
\bigg |\frac{\partial K_i}{\partial y^1}(0) \bigg |\int  _{B_{\sigma}} v_i^{p_i+1} dy & \leq & C \bigg| \int _{B_{\sigma}} \lambda _i^{\frac{2n}{n-2}}v_i^{p_i+1}\frac{\partial K_i}{\partial y^1}(0) dy \bigg|  \\
& \leq & C \left (r_iv_i(0)^{-\frac{2}{n-2} }+\left (\delta_i r_i
+ v_i(0)^{2 t _i} \right ) \right )\\
&  \leq & C \left (r_iv_i(0)^{-\frac{2}{n-2}} + r_i v_i(0)^{2t_i}
+ v_i(0)^{2t_i} \right ) \hspace{.2in} \text{ (by inequality (\ref{eq:deltaiprelim}))}\\
& \leq & C \left (r_iv_i(0)^{-\frac{2}{n-2}}+ v_i(0)^{2t_i} \right ).
\end{eqnarray*} 
Then by (\ref{eq:lower-int-v^p+1}) 
\begin{equation}
\label{eq:gradKprelim}
\bigg |\frac{\partial K_i}{\partial y^1}(0) \bigg | \leq C
\left (r_iv_i(0)^{-\frac{2}{n-2}} + v_i(0)^{2 t _i}
\right ).
\end{equation}

The same estimate holds for $\big |\frac{\partial K_i}{\partial
y^j}(0) \big |, \,\, j=2,...,n $ as well, since we can also choose $X=\frac{\partial}{\partial y^j}$ in the above calculation.

\subsection{Location of the Blow-up}
\label{subsection:location}
  
Choose a point $\bar{y}$ with $|\bar{y}|=1$.
It is proved in Section 6 of \cite{YY3} that $\frac{v_i}{v_i(\bar{y})}$ converges in
$C^2$-norm to a function $h$ on any compact subset of $\mathbf{R}^n \setminus \{0\}$, and $h=\frac{1}{2} + \frac{1}{2}|y|^{2-n}$.

Recall that we chose the coordinate systems $z=(z^1,...,z^n)$ and $y=\frac{z}{r_i}$ to be centered at each $x_i \in M$, thus $\nabla K_i(0)=r_i\nabla K(x_i)$.  Here we write $\nabla K(x_i)$ instead of $\nabla K(0)$ to emphasize the fact that $\nabla K$ is evaluated at different point $x_i$ as $i \to \infty$.  We claim that this blow-up must occur at a critical point of $K$, i.e.,

\begin{prop}
\label {prop:gradK=0}
$\nabla K(x_0)=\displaystyle \lim_{i \to \infty} \nabla K(x_i)=0$.
\end{prop}

\pf 
Suppose this is not true, then there exists some $j \in \{1,...,n\}$, such that$\big | \frac{\partial K}{\partial z^j}(x_i) \big | \geq \varepsilon$ for a constant $\varepsilon$ independent of $i$.  Without loss of generality we assume $j=1$.  Then from inequality (\ref{eq:gradKprelim}) we know that 
$\varepsilon r_i \leq C
\left (r_iv_i(0)^{-\frac{2}{n-2}} + v_i(0)^{2 t _i}\right ).  $
Therefore 
\begin{equation}
\label{eq:location-ri}
r_i \leq Cv_i(0)^{2 t _i}
\end{equation}
when $v_i(0)^{-\frac{2}{n-2}}$ is sufficiently small.

Once more we look at the Pohozaev identity (\ref{eq:pohozaev-X(R)}) with $X=\displaystyle \sum _{j}y^j \frac{\partial}{\partial y^j}$.  We divide both sides of it by $v_i^2(\bar{y})$ so it becomes 
\begin{equation}
\label{eq:pohozaev-case1}
\frac{n-2}{2n} \frac{1}{v_i^2(\bar{y})}\int _{B_{\sigma}} X(R_i) dv _{g_i} =\frac{1}{v_i^2(\bar{y})} \int _{\partial
B_{\sigma}} T_i(X, \nu _i) d \Sigma _i
\end{equation}

Its right hand side is 

\begin{eqnarray}
\label{eq:case1-Tterm}
 & & \frac{1}{v_i^2(\bar{y})}\int _{\partial
B_{\sigma}} T_i(X, \nu _i) d \Sigma _i \nonumber \\
& = & \frac{1}{v_i^2(\bar{y})}\int _{\partial
B_{\sigma}}\left ( \Ric (g_i)-n^{-1}R(g_i)g_i \right )(X, \nu _i) d \Sigma _i \nonumber \\
& = & \frac{1}{v_i^2(\bar{y})}\int _{\partial B_{\sigma}} 
\bigg [ \Ric 
  \left ( 
    \left (\lambda _i v_i \right )^{\frac{4}{n-2}}dy \otimes dy \right ) \nonumber \\
& & - n^{-1} R 
  \left( \left(\lambda _i v_i \right)^{\frac{4}{n-2}}dy \otimes dy \right) 
  \left( \lambda _i v_i \right )^{\frac{4}{n-2}}dy \otimes dy
\bigg ](X, \nu_0) (\lambda _i v_i)^2 d \Sigma _{\sigma} \nonumber \\
& = &  \int _{\partial B_{\sigma}} 
\left( \frac{\lambda _i v_i}{v_i(\bar{y})}\right)^2 
\Bigg [ \Ric
  \left ( 
     \left(\frac{\lambda _i v_i}{v_i(\bar{y})}\right)^{\frac{4}{n-2}} dy \otimes dy
  \right) \\
  & & - n^{-1} R 
  \left ( 
    \left(\frac{\lambda _i
v_i}{v_i(\bar{y})}\right)^{\frac{4}{n-2}}dy \otimes dy
  \right) 
  \left( \frac{\lambda _i v_i}{v_i(\bar{y})} \right ) 
  ^{\frac{4}{n-2}}dy \otimes dy
\Bigg ] (X, \nu _0)    d \Sigma_{\sigma} \nonumber 
\end{eqnarray} 
 where $\nu _0 =\sigma ^{-1} \displaystyle \sum _{j}y^j \frac{\partial}{\partial y^j}$ is the unit outer normal on
$\partial B_{\sigma}$ with respect to 
the Euclidean metric $dy \otimes dy$. 

\noindent
When $i \rightarrow \infty$, for $|y|=\sigma$, $\lambda _i(y) =
\lambda (r _i y) \rightarrow \lambda (x_0)$. Thus when $i$ goes to $\infty$, up to a constant
(\ref{eq:case1-Tterm}) converges to 
\begin{eqnarray}
\label{eq:case1-Tterm-followup}
& & \int _{\partial
B_{\sigma}} h^2  \bigg (\Ric \left (h^{\frac{4}{n-2}}dy \otimes dy \right )  - n^{-1}
R \left (h^{\frac{4}{n-2}} dy \otimes dy \right )h^{\frac{4}{n-2}}  dy \otimes dy \bigg )(X, \nu _0)    d \Sigma _{\sigma}       \nonumber \\
& = &  \int _{\partial
B_{\sigma}} h^2 \cdot (n-2) h^{\frac{2}{n-2}} \left [\Hess \left
(h^{-\frac{2}{n-2}} \right) (X, \nu _0) -\frac{1}{n}\Delta \left (h^{-\frac{2}{n-2}} \right )<X, \nu _0>\right ]   \,\,  d \Sigma _{\sigma}
\nonumber \\
& = & (n-2)\sigma^{-1} \int _{\partial
B_{\sigma}} h^{\frac{2(n-1)}{n-2}} \cdot \left [ \Hess \left
(h^{-\frac{2}{n-2}} \right )(X,X)  -\frac{1}{n}\Delta \left (h^{-\frac{2}{n-2}} \right )\sigma ^2  \right ] \,\,  d \Sigma _{\sigma}
\end{eqnarray}

\noindent
We know that 
\begin{equation*}
h^{-\frac{2}{n-2}} \,\, = \,\, \left( \frac{1}{2} (1+|y|^{2-n}) \right )^{-\frac{2}{n-2}} \,\, = \,\,  2^{\frac{2}{n-2}}|y|^2 - \frac{2^{\frac{n}{n-2}}}{n-2}|y|^{n} +
O \left (|y|^{2(n-1)} \right ),
\end{equation*}

\noindent
and by direct computation 
$$\Hess \left( 2^{\frac{2}{n-2}}|y|^2 - \frac {2^{\frac{n}{n-2}}}{n-2}
|y|^{n}  \right )(X,X) -
\frac{1}{n}\Delta \left ( 2^{\frac{2}{n-2}}|y|^2 - \frac {2^{\frac{n}{n-2}}}{n-2}
|y|^{n} \right ) \sigma ^2 = -2^{\frac{n}{n-2}}(n-1)  \sigma ^n .$$

\noindent
Therefore $$ \Hess \left (h^{-\frac{2}{n-2}} \right )(X,X)
-\frac{1}{n}\Delta \left (h^{-\frac{2}{n-2}} \right )\sigma ^2 =  -2^{\frac{n}{n-2}}(n-1) \sigma
^n + O \left (\sigma ^{2(n-1)} \right )  . $$    
Also we know
 \begin{eqnarray*}
 h^{\frac{2(n-1)}{n-2}}
& = & \left (\frac
{1}{2} \right )^{\frac{2(n-1)}{n-2}}|y|^{-2(n-1)}\left (1+ O(|y|^{n-2}) \right ).
\end{eqnarray*}

\noindent
Thus we can conclude that (\ref{eq:case1-Tterm-followup}) is equal to 
\begin{eqnarray}
\label{eq:case1-Tterm-<0}
 &  & -\frac{1}{2}(n-1)(n-2) \sigma^{-1} \int _{\partial B_{\sigma} }
\left (|y|^{-2(n-1)}+O(|y|^{-n})\right )  \left ( |y|^{n}+O(|y|^{2(n-1)})
\right )\sigma^{n-1}d\Sigma
_1 \nonumber \\
& = &  -\frac{1}{2}(n-1)(n-2) + O(\sigma ^{n-2}) \nonumber 
\end{eqnarray}
Therefore the limit of the right hand side of (\ref{eq:pohozaev-case1}) is strictly less than $0$
when we choose $\sigma$ to be sufficiently small. 

On the other hand, the left hand side of (\ref{eq:pohozaev-case1}) is 
$$
\frac{n-2}{2n}c(n)^{-1} \frac{1}{v^2_i(\bar{y})}\int _{B_{\sigma}}
X(K_iv_i^{-\delta _i})(\lambda _i v_i)^{\frac{2n}{n-2}}dy.
$$
We write
\begin{eqnarray}
\label{eq:case1-X(R)term}
& & \frac{1}{v^2_i(\bar{y})}\int _{B_{\sigma}}
X(K_iv_i^{-\delta _i })(\lambda _i
v_i)^{\frac{2n}{n-2}}dy \nonumber \\
 & =
& \frac{1}{v^2_i(\bar{y})}\int _{B_{\sigma}}
X(K_i)v_i^{p_i +1 }\lambda _i ^{\frac{2n}{n-2}}dy
-\frac{\delta _i}{v^2_i(\bar{y})}
  \int _{B_{\sigma}} K_i\lambda _i ^{\frac{2n}{n-2}}v_i^{p_i}X(v_i)
dy .
\end{eqnarray}
The second term of ({\ref {eq:case1-X(R)term})

\begin{eqnarray*}
& = & -\frac{\delta _i}{p_i+1}\frac{1}{v^2_i(\bar{y})}
   \int _{B_{\sigma}} K_i\lambda _i ^{\frac{2n}{n-2}}X(v_i^{p_i+1})
dy \\
& = &  -\frac{\delta _i}{p_i+1}\frac{1}{v^2_i(\bar{y})} \int
_{B_{\sigma}} \bigg [\diver \left (K_i \lambda _i ^{\frac{2n}{n-2}}v_i^{p_i+1} X \right ) -K_i \lambda _i ^{\frac{2n}{n-2}}v_i^{p_i+1} \diver X   \\
& &  -\lambda _i
^{\frac{2n}{n-2}}v_i^{p_i+1}X(K_i) - K_i
v_i^{p_i+1} X (\lambda _i ^{\frac{2n}{n-2}}) \bigg ]
dy \\
& = &  -\frac{\delta _i}{p_i+1}\frac{\sigma}{v^2_i(\bar{y})} \int
_{\partial B_{\sigma}} K_i \lambda _i
^{\frac{2n}{n-2}}v_i^{p_i+1} d \Sigma _{\sigma} \\
& &  + \frac{\delta
_i}{p_i+1}\frac{1}{v^2_i(\bar{y})} \int _{B_{\sigma}} K_i \lambda _i
^{\frac{2n}{n-2}}v_i^{p_i+1} \left (n + X(\ln K_i)+ \frac{2n}{n-2}
X(\ln \lambda _i) \right ) dy
\end{eqnarray*}

\noindent
On $\partial B_{\sigma}$, $\frac{v_i}{v_i(\bar{y})}
\rightarrow h(\sigma)$ and $v_i \rightarrow 0$ uniformly, so
$$
\frac{1}{v^2_i(\bar{y})} \int
_{\partial B_{\sigma}} K_i \lambda _i
^{\frac{2n}{n-2}}v_i^{p_i+1} d \Sigma _{\sigma}   =   \int
_{\partial B_{\sigma}} K_i \lambda _i
^{\frac{2n}{n-2}} \left ( \frac {v_i}{v_i(\bar{y})} \right )^2  v_i^{p_i-1} d \Sigma
_{\sigma}  \rightarrow  0 .
$$

\noindent
Since $X=r\frac{\partial}{\partial r}$ and $\big | \frac{\partial}{\partial
r}(\ln K_i) \big | $, $\big | \frac{\partial}{\partial
r}(\ln \lambda _i) \big |$ are uniformly bounded, we can choose $\sigma$ to be
small (independent of $i$) to make $n + X(\ln K_i )+ \frac{2n}{n-2}
X(\ln \lambda _i)>0$.  Thus when $i \rightarrow \infty$, the limit of the second term of (\ref{eq:case1-X(R)term}) is greater than or equal
to $0$.  

Next we will show that the limit of the first
term of (\ref{eq:case1-X(R)term}) is $0$, or equivalently,
\begin{equation}
\label{eq:key-X(Ki)}
\lim _{i \rightarrow \infty} v^2_i(0)\int _{B_{\sigma}}
X(K_i)v_i^{p_i+1}\lambda _i^{\frac{2n}{n-2}}dy
=0,
\end{equation}
since $v_i(\bar{y})\geq Cv_i(0)^{-1}$ by Proposition \ref{prop:simpleestimates}.
This then will end the proof because it implies that the limit of the left hand side of (\ref{eq:pohozaev-case1}) is greater than
or equal to $0$, contradicting the sign of the right hand side.

\noindent
Note that 
\begin{eqnarray*}
X(K_i)(y) 
 & = & \left (\sum _{j} y^j \frac{\partial
K_i}{\partial y^j} \right )(y) \\
& = & \left (\sum _{j} y^j \frac{\partial K_i}{\partial y^j} \right )(0) + \sum _{k} 
\frac{\partial}{\partial y^k} \left (\sum _{j}y^j \frac{\partial K_i}{\partial y^j} \right ) (\varsigma)y^k \hspace{.2in} \text{for some } |\varsigma| \leq |y| \\
& = & \sum _{j} \frac{\partial K_i}{\partial y^j}(\varsigma) y^j + \sum _{j,k} \frac{\partial^2 K_i}{\partial y^k \partial y^j}(\varsigma) \varsigma ^j y^k  
\end{eqnarray*}

\noindent
Therefore 
\begin{eqnarray*}
& & v^2_i(0) \bigg | \int _{B_{\sigma}}
X(K_i)v_i^{p_i+1}\lambda _i^{\frac{2n}{n-2}} dy\bigg |\\
& \leq &   v_i^2(0) \int _{B_{\sigma}} \sum _{j} \bigg | \frac{\partial K_i}{\partial y^j}(\varsigma) \bigg | |y| v_i^{p_i+1}\lambda _i^{\frac{2n}{n-2}}dy +  v_i^2(0) \int _{B_{\sigma}} \sum _{j,k} \bigg | \frac{\partial ^2 K_i}{\partial y^j \partial y^k}(\varsigma) \bigg | |y|^2 v_i^{p_i+1}\lambda _i^{\frac{2n}{n-2}}dy \\
& \leq &  Cv_i^2(0) r_i \int _{B_{\sigma}} |y| v_i^{p_i+1}dy + C v_i^2(0) r_i ^2\int _{B_{\sigma}} |y|^2 v_i^{p_i+1}dy \\
& \leq & Cv_i^2(0)r_i\cdot v_i(0)^{-\frac{2}{n-2}+\frac{n-1}{2}\delta _i}+Cv_i^2(0)r_i^2\cdot v_i(0)^{-\frac{4}{n-2}+\frac{n}{2}\delta _i} \hspace {.2in} \text{ (by Lemma } \ref{lemma:intv^p+1y^k})\\
& \leq & Cv_i(0)^{2+2t_i-\frac{2}{n-2}}+Cv_i(0)^{2+4t_i-\frac{4}{n-2}} \hspace {.2in} \text{ (by Proposition } \ref {prop:v^delta} \text{ and Inequality } (\ref {eq:location-ri})) 
\end{eqnarray*}

\noindent
By the definition of $t_i$,
\begin{equation*}
\lim _{i \to \infty}t_i \,\, = \,\, \lim _{i \to \infty}\left( 1-\frac{(p_i-1)l_i}{2}\right) \,\, = \,\, 1-\frac{2}{n-2}\lim _{i \to \infty} l_i \,\, < \,\, 1- \frac{2}{n-2}\cdot \frac{2n-5}{2} = \frac{3-n}{n-2}
\end{equation*}

\noindent
Thus
\begin{equation}
\label {eq:limit2+2ti-2/n-2}
\lim _{i \to \infty} \left ( 2+2t_i-\frac{2}{n-2} \right ) \,\, < \,\,  2 +2 \cdot
 \frac{3-n}{n-2}  - \frac{2}{n-2} \,\, =  \,\, 0 
\end{equation}
and
\begin{equation*}
\lim _{i \to \infty} \left ( 2+4t_i-\frac{4}{n-2} \right ) \,\, < \,\,  2 +4 \cdot
 \frac{3-n}{n-2}  - \frac{4}{n-2} \,\, = \,\,\frac{4-2n}{n-2} \,\, <  \,\, 0 .
\end{equation*}

\noindent
Since these are all strict inequalities, we know that 
\begin{equation*}
\lim _{i \to \infty} \left( Cv_i(0)^{2+2t_i-\frac{2}{n-2}}+Cv_i(0)^{2+4t_i-\frac{4}{n-2}} \right )=0
\end{equation*}
and consequently 

$$\lim_{i \to \infty}  v^2_i(0) \bigg | \int _{B_{\sigma}}
X(K_i)v_i^{p_i+1}\lambda _i^{\frac{2n}{n-2}} dy\bigg |=0.$$

\stop

\subsection{Refined Estimates for $\delta_i$ and $|\nabla K_i|$}
\label{subsection:refined}

Now because $x_0=\displaystyle \lim_{i \to \infty} x_i$ is a critical point of the function $K$, which satisfies the flatness condition $(*)$, we have 
$|\nabla ^p K(x_i) |\leq C_0 | \nabla K(x_i) |^{\frac{n-2-p}{n-3}}$ when $2 \leq p \leq n-3$.  When $p=2$, since $g=\lambda ^{\frac{4}{n-2}} dz^2$, this implies

$$\bigg | \nabla ^2 K \left (  \frac {\partial}{\partial z^{l_1}},  \frac {\partial}{\partial z^{l_2}} \right ) (x_i)   \bigg | \,\, = \,\, \bigg | \frac{\partial ^2K}{\partial z^{l_1}\partial z^{l_2}} (x_i)-\Gamma _{l_1l_2}^l(x_i)\frac{\partial K}{\partial z ^l}(x_i) \bigg | \,\, \leq \,\, C |\nabla K(x_i)|^{\frac{n-4}{n-3}},  $$ where $l_1, l_2, l =1,2,...,n.$  Therefore 
$$ \bigg | \frac{\partial ^2K}{\partial z^{l_1}\partial z^{l_2}} (x_i) \bigg | \,\, \leq \,\, C |\nabla K (x_i)|+C |\nabla K (x_i)|^{\frac{n-4}{n-3}} \,\, \leq C |\nabla K (x_i)|^{\frac{n-4}{n-3}},$$ since $|\nabla K (x_i)| <1$ for sufficiently large $i$.  That is,
$\big | \frac{\partial ^{\alpha} K}{\partial z^{\alpha}}(x_i) \big | \leq C | \nabla K (x_i) |^{\frac{n-2-|\alpha|}{n-3}}$ for $|\alpha|=2.$  Here we have used the notations that  

$$\alpha=(\alpha _1, \alpha _2,..., \alpha _n) \text{ with each } \alpha _i \geq 0, \hspace{.2in} |\alpha|=\alpha_1+\alpha_2+...+\alpha _n,$$ and $$ \frac{\partial ^{\alpha} K}{\partial z^{\alpha}}  = \frac{\partial^ {\alpha_1}\partial^ {\alpha_2}...\partial^ {\alpha_n} K}{(\partial {z^1})^{\alpha _1}(\partial {z^2})^{\alpha _2}...(\partial {z^n})^{\alpha _n} }.$$

\noindent
Generally, when $2 \leq p < q \leq n-3$, we have $|\nabla K (x_i)|^{\frac{n-2-p}{n-3}} < |\nabla K (x_i)|^{\frac{n-2-q}{n-3}}$, so by similar computations we have 
$$\bigg | \frac{\partial ^{\alpha} K}{\partial z^{\alpha}}(x_i) \bigg | \leq C | \nabla K(x_i) |^{\frac{n-2-|\alpha|}{n-3}}  \hspace{.2in} \text{for } \hspace{.2in} 2\leq |\alpha|\leq n-3 .$$  

\noindent
Then since $K_i(y)=K(r_iy)$, $\big | \frac{\partial ^{\alpha} K_i}{\partial y^{\alpha}}(0)\big|=r_i^{|\alpha|}\big | \frac{\partial ^{\alpha} K}{\partial z^{\alpha}}(x_i)\big|$ and $ | \nabla K_i(0) |=r_i | \nabla K(x_i) |$.
Thus 
\begin{eqnarray}
\label{eq:higherderiofK}
\bigg | \frac{\partial ^{\alpha} K_i}{\partial y^{\alpha}}(0)\bigg| & \leq & r_i^{|\alpha|}C  | \nabla K(x_i) |^{\frac{n-2-|\alpha|}{n-3}}\nonumber \\
& = & Cr_i^{\frac{(|\alpha|-1)(n-2)}{n-3}} | \nabla K_i(0) |^{\frac{n-2-|\alpha|}{n-3}} \nonumber \\
& < & C r_i  | \nabla K_i(0) |^{\frac{n-2-|\alpha|}{n-3}} 
\end{eqnarray}
where the last step comes from the fact that $\frac{(|\alpha|-1)(n-2)}{n-3}>1$ and $r_i<1$.  With this flatness condition on $K_i$, we can refine the estimates for $\delta_i$ and $|\nabla K_i|$} as follows.  

\noindent
Inequality (\ref {eq:deltacoarse}) gives
$$\delta _i \leq C \left ( v_i(0)^{2t _i} + \bigg | \int  _{B_{\sigma}}   \frac {\partial K_i}{\partial r}|y| v_i^{p_i+1}\lambda _i^{\frac{2n}{n-2}}dy \bigg | \right )=C\left ( v_i(0)^{2t _i} + \bigg | \int  _{B_{\sigma}}  r \frac {\partial K_i}{\partial r} v_i^{p_i+1}\lambda _i^{\frac{2n}{n-2}}dy \bigg| \right ) .$$  

\noindent
We write $ r \frac {\partial K_i}{\partial r}=\displaystyle \sum_{j}y^j\frac {\partial K_i}{\partial y^j}$.  For each $j=1,...,n$,
\begin{eqnarray*}
\frac {\partial K_i}{\partial y^j}(y) & = & \frac {\partial K_i}{\partial y^j}(0) + \sum _{|\beta|=1} \frac{\partial ^{\beta}}{\partial y^{\beta}}\frac {\partial K_i}{\partial y^j}(0)y^{\beta}+\frac{1}{2!}\sum _{|\beta|=2} \frac{\partial ^{\beta}}{\partial y^{\beta}}\frac {\partial K_i}{\partial y^j}(0)y^{\beta}+\cdot \cdot \cdot \\
& & + \frac{1}{(n-4)!}\sum _{|\beta|=n-4} \frac{\partial ^{\beta}}{\partial y^{\beta}}\frac {\partial K_i}{\partial y^j}(0)y^{\beta}+\frac{1}{(n-3)!}\sum _{|\beta|=n-3} \frac{\partial ^{\beta}}{\partial y^{\beta}}\frac {\partial K_i}{\partial y^j}(\varsigma)y^{\beta}
\end{eqnarray*}

\noindent
where $|\varsigma|\leq |y|$, and $y^{\beta}=y_1^{\beta _1}y_2^{\beta _2}\cdot \cdot \cdot y_n^{\beta_n}$ for $\beta=(\beta_1, \beta_2,...,\beta_n)$.  Therefore
\begin{eqnarray*}
& &  \int  _{B_{\sigma}}\bigg |  r \frac {\partial K_i}{\partial r}\bigg | v_i^{p_i+1}\lambda _i^{\frac{2n}{n-2}}dy  \\
 & \leq & C \Bigg ( \int _{B_{\sigma}} \bigg|  \frac {\partial K_i}{\partial y^j}(0)\bigg| |y|v_i^{p_i+1}dy + \sum _{|\beta|=1}^{n-4} \int _{B_{\sigma}} \bigg | \frac{\partial ^{\beta}}{\partial y^{\beta}}\frac {\partial K_i}{\partial y^j}(0)\bigg| |y|^{|\beta|+1} v_i^{p_i+1}dy\\
& & + \sum _{|\beta|=n-3} \int _{B_{\sigma}} \bigg | \frac{\partial ^{\beta}}{\partial y^{\beta}}\frac {\partial K_i}{\partial y^j}(\varsigma)\bigg| |y|^{n-2} v_i^{p_i+1}dy \Bigg ).
\end{eqnarray*}

\noindent
By Lemma \ref{lemma:intv^p+1y^k} and Proposition \ref {prop:v^delta}, the first term $$\int _{B_{\sigma}} \bigg |  \frac {\partial K_i}{\partial y^j}(0)\bigg | |y|v_i^{p_i+1}dy \leq C   |\nabla K_i(0)  | v_i(0)^{-\frac{2}{n-2}},$$ and the last term $$\sum _{|\beta|=n-3} \int _{B_{\sigma}} \bigg | \frac{\partial ^{\beta}}{\partial y^{\beta}}\frac {\partial K_i}{\partial y^j}(\varsigma)\bigg| |y|^{n-2} v_i^{p_i+1}dy \leq Cr_i^{n-2}v_i(0)^{-2}. $$

\noindent
In addition, by (\ref {eq:higherderiofK}), for any $1 \leq |\beta| \leq n-4$, 

\begin{eqnarray}
\label {eq:int-gradK-y^-v^p+1}
& & \int _{B_{\sigma}} \bigg | \frac{\partial ^{\beta}}{\partial y^{\beta}}\frac {\partial K_i}{\partial y^j}(0)\bigg| |y|^{|\beta|+1} v_i^{p_i+1}dy \nonumber \\
 & \leq & C r_i\int _{B_{\sigma}}  | \nabla K_i(0) |^{\frac{n-2-(|\beta|+1)}{n-3}}   |y|^{|\beta|+1} v_i^{p_i+1}dy \nonumber \\
& = & Cr_i \int _{B_{\sigma}} | \nabla K_i(0) |^{\frac{n-3-|\beta|}{n-3}}|y|^{|\beta|}\cdot|y|v_i^{p_i+1} dy \nonumber \\
& \leq & Cr_i \int _{B_{\sigma}}\left ( | \nabla K_i(0) |^{\frac{n-3-|\beta|}{n-3}\cdot \frac{n-3}{n-3-|\beta|}}+ |y|^{|\beta|\cdot \frac{n-3}{|\beta|} } \right )\cdot|y|v_i^{p_i+1} dy \nonumber \\
& & \text{ (by Young's Inequality)} \nonumber \\
& = & Cr_i \left( \int _{B_{\sigma}}  | \nabla K_i(0) | \cdot |y|v_i^{p_i+1} dy + \int _{B_{\sigma}} |y|^{n-2}v_i^{p_i+1} dy  \right) \nonumber \\
& \leq &  Cr_i | \nabla K_i(0) | v_i(0)^{-\frac{2}{n-2}}+ Cr_iv_i(0)^{-2}. \nonumber
\end{eqnarray}
Thus 
\begin{eqnarray}
\label {eq:int-rdK/dr-v^p+1}
& & \int  _{B_{\sigma}}  \bigg | r \frac {\partial K_i}{\partial r}\bigg | v_i^{p_i+1} \lambda _i^{\frac{2n}{n-2}}dy \nonumber  \\
& \leq &   C  | \nabla K_i (0)| v_i(0)^{-\frac{2}{n-2}}+ \left ( Cr_i | \nabla K_i(0) | v_i(0)^{-\frac{2}{n-2}}+ Cr_iv_i(0)^{-2} \right  )  + Cr_i^{n-2}v_i(0)^{-2}\nonumber \\
& \leq & C |\nabla K_i (0)| v_i(0)^{-\frac{2}{n-2}}+ Cr_iv_i(0)^{-2}.
\end{eqnarray}

\noindent
Plugging this back into (\ref {eq:deltacoarse}) we now have a refined estimate
\begin{equation}
\label {eq:delta-refined}
\delta _i \leq C \left ( v_i(0)^{2t_i} + |\nabla K_i (0)| v_i(0)^{-\frac{2}{n-2}}+ r_iv_i(0)^{-2}  \right ).
\end{equation}

This will enable us to also refine the estimate for $|\nabla K _i(0)|$.

\noindent
Inequality (\ref {eq:gradKicoarse}) gives $\big| \int _{B_{\sigma}} \lambda _i^{\frac{2n}{n-2}} v_i^{p_i+1}\frac{\partial K_i}{\partial
y^1} dy \big| \leq C \left ( \delta _i r_i + v_i(0)^{2 t _i} \right ).$

\noindent
Again we write
\begin{eqnarray*}
\frac {\partial K_i}{\partial y^1}(y) & = & \frac {\partial K_i}{\partial y^1}(0) + \sum _{|\beta|=1} \frac{\partial ^{\beta}}{\partial y^{\beta}}\frac {\partial K_i}{\partial y^1}(0)y^{\beta}+\frac{1}{2!}\sum _{|\beta|=2} \frac{\partial ^{\beta}}{\partial y^{\beta}}\frac {\partial K_i}{\partial y^1}(0)y^{\beta}+\cdot \cdot \cdot \\
& & + \frac{1}{(n-4)!}\sum _{|\beta|=n-4} \frac{\partial ^{\beta}}{\partial y^{\beta}}\frac {\partial K_i}{\partial y^1}(0)y^{\beta}+\frac{1}{(n-3)!}\sum _{|\beta|=n-3} \frac{\partial ^{\beta}}{\partial y^{\beta}}\frac {\partial K_i}{\partial y^1}(\varsigma)y^{\beta}.
\end{eqnarray*}

\noindent
Therefore we have 
\begin{eqnarray*}
& & \int _{B_{\sigma}}\lambda _i^{\frac{2n}{n-2}} v_i^{p_i+1} \bigg | \frac{\partial K_i}{\partial
y^1}(0) \bigg | dy \\
& \leq & \bigg | \int _{B_{\sigma}} \lambda _i^{\frac{2n}{n-2}}v_i^{p_i+1}\frac{\partial K_i}{\partial y^1} dy \bigg |  + C\sum _{|\beta|=1}^{n-4} \int _{B_{\sigma}} \bigg | \frac{\partial ^{\beta}}{\partial y^{\beta}}\frac {\partial K_i}{\partial y^1}(0)\bigg| |y|^{|\beta|} v_i^{p_i+1}dy\\
& & +C \sum _{|\beta|=n-3} \int _{B_{\sigma}} \bigg | \frac{\partial ^{\beta}}{\partial y^{\beta}}\frac {\partial K_i}{\partial y^1}(\varsigma)\bigg| |y|^{n-3} v_i^{p_i+1}dy \\
 & \leq &  C \left ( \delta _i r_i + v_i(0)^{2 t _i} \right ) + C\sum _{|\beta|=1}^{n-4} \int _{B_{\sigma}} \bigg | \frac{\partial ^{\beta}}{\partial y^{\beta}}\frac {\partial K_i}{\partial y^1}(0)\bigg| |y|^{|\beta|} v_i^{p_i+1}dy\\
& & +C \sum _{|\beta|=n-3} \int _{B_{\sigma}} \bigg | \frac{\partial ^{\beta}}{\partial y^{\beta}}\frac {\partial K_i}{\partial y^1}(\varsigma)\bigg| |y|^{n-3} v_i^{p_i+1}dy .
\end{eqnarray*}
By (\ref {eq:lower-int-v^p+1}) this implies 
\begin{eqnarray*}
\bigg | \frac{\partial K_i}{\partial
y^1}(0) \bigg | &  \leq &   C \left ( \delta _i r_i + v_i(0)^{2 t _i} \right ) + C\sum _{|\beta|=1}^{n-4} \int _{B_{\sigma}} \bigg | \frac{\partial ^{\beta}}{\partial y^{\beta}}\frac {\partial K_i}{\partial y^1}(0)\bigg| |y|^{|\beta|} v_i^{p_i+1}dy \\
& & +C \sum _{|\beta|=n-3} \int _{B_{\sigma}} \bigg | \frac{\partial ^{\beta}}{\partial y^{\beta}}\frac {\partial K_i}{\partial y^1}(\varsigma)\bigg| |y|^{n-3} v_i^{p_i+1}dy .
\end{eqnarray*}

\noindent
By Lemma \ref{lemma:intv^p+1y^k}, Proposition \ref {prop:v^delta}, (\ref {eq:higherderiofK}), and Young's Inequality, when $1 \leq |\beta| \leq n-4$,
\begin{eqnarray*}
& & \int _{B_{\sigma}} \bigg | \frac{\partial ^{\beta}}{\partial y^{\beta}}\frac {\partial K_i}{\partial y^1}(0)\bigg| |y|^{|\beta|} v_i^{p_i+1}dy  \\
 & \leq & C r_i\int _{B_{\sigma}}  | \nabla K_i(0) |^{\frac{n-2-(|\beta|+1)}{n-3}}   |y|^{|\beta|} v_i^{p_i+1}dy \nonumber \\
& = & Cr_i \int _{B_{\sigma}} | \nabla K_i(0) |^{\frac{n-3-|\beta|}{n-3}}|y|^{|\beta|}v_i^{p_i+1} dy  \\
& \leq & Cr_i \int _{B_{\sigma}}\left ( | \nabla K_i(0) |^{\frac{n-3-|\beta|}{n-3}\cdot \frac{n-3}{n-3-|\beta|}}+ |y|^{|\beta|\cdot \frac{n-3}{|\beta|} } \right )v_i^{p_i+1} dy  \\
& = & Cr_i \left( \int _{B_{\sigma}}  | \nabla K_i(0) |  v_i^{p_i+1} dy + \int _{B_{\sigma}} |y|^{n-3}v_i^{p_i+1} dy  \right) \\
& \leq &  Cr_i | \nabla K_i(0) | + Cr_iv_i(0)^{-\frac{2(n-3)}{n-2}}.   
\end{eqnarray*}
Furthermore, $$\sum _{|\beta|=n-3} \int _{B_{\sigma}} \bigg | \frac{\partial ^{\beta}}{\partial y^{\beta}}\frac {\partial K_i}{\partial y^1}(\varsigma)\bigg| |y|^{n-3} v_i^{p_i+1}dy \,\,  \leq \,\, Cr_i^{n-2}\int _{B_{\sigma}}  |y|^{n-3} v_i^{p_i+1}dy \,\, \leq \,\,  Cr_i^{n-2}v_i(0)^{-\frac{2(n-3)}{n-2}}. $$

\noindent
Therefore
\begin{eqnarray*}
\bigg | \frac{\partial K_i}{\partial y^1}(0) \bigg | &  \leq &   
C \left (\delta _i r_i+v_i(0)^{2t_i} \right )
+ \left (  Cr_i | \nabla K_i(0) | + Cr_iv_i(0)^{-\frac{2(n-3)}{n-2}} \right )  +  Cr_i^{n-2}v_i(0)^{-\frac{2(n-3)}{n-2}} \\
& \leq & C \delta _i r_i + C v_i(0)^{2t _i} +   Cr_i | \nabla K_i(0) | + Cr_iv_i(0)^{-\frac{2(n-3)}{n-2}}.  
\end{eqnarray*}
The same estimate also holds for $\big | \frac{\partial K_i}{\partial y^j}(0) \big |$, where $j=2,..., n$, so we know
\begin{eqnarray*}
|\nabla K_i(0)  | & \leq & C \delta _i r_i + C v_i(0)^{2t _i} +   Cr_i | \nabla K_i(0) | + Cr_iv_i(0)^{-\frac{2(n-3)}{n-2}}  \\
&  \leq &   C \left ( v_i(0)^{2t_i} + |\nabla K_i (0)| v_i(0)^{-\frac{2}{n-2}}+ r_iv_i(0)^{-2}  \right ) r_i + C v_i(0)^{2t _i}\\
& &  +   Cr_i | \nabla K_i(0) | + Cr_iv_i(0)^{-\frac{2(n-3)}{n-2}} \hspace{.4in} \text{( by (\ref {eq:delta-refined}) )}.
\end{eqnarray*}
When $i$ is large enough, all the terms involving  $|\nabla K_i (0)|$ can be absorbed into the left hand side of this inequality, therefore we get a refined estimate 
\begin{eqnarray}
\label {eq:gradK-refined}
 |\nabla K_i (0)| & \leq & Cr_iv_i(0)^{2t_i}+Cr_i^2v_i(0)^{-2}+C v_i(0)^{2t _i} + Cr_iv_i(0)^{-\frac{2(n-3)}{n-2}} \nonumber \\
& \leq & Cr_i^2v_i(0)^{-2}+C v_i(0)^{2t _i} + Cr_iv_i(0)^{-\frac{2(n-3)}{n-2}} .
\end{eqnarray}

Finally, we are going to prove that (\ref {eq:key-X(Ki)}) holds.  As in the proof of Proposition \ref {prop:gradK=0}, this will give the desired contradiction by comparing the signs of both sides of (\ref {eq:pohozaev-case1}), which rules out case I.

\noindent
We know
\begin{eqnarray*}
& &  v^2_i(0)\int _{B_{\sigma}} | X(K_i)| v_i^{p_i+1}\lambda _i^{\frac{2n}{n-2}}dy \\
& = & v^2_i(0)\int _{B_{\sigma}} \bigg | r\frac{\partial K_i}{\partial r}  \bigg| v_i^{p_i+1}\lambda _i^{\frac{2n}{n-2}}dy \\
& \leq & C v^2_i(0) \left ( |\nabla K_i (0)| v_i(0)^{-\frac{2}{n-2}}+r_i v_i(0)^{-2} \right ) \hspace{.2in} \text {\big ( by (\ref {eq:int-rdK/dr-v^p+1}) \big )} \\
& \leq &  C v^2_i(0) \Bigg ( \left (  r_i^2v_i(0)^{-2}+ v_i(0)^{2t _i} + r_iv_i(0)^{-\frac{2(n-3)}{n-2}} \right ) v_i(0)^{-\frac{2}{n-2}}+r_i v_i(0)^{-2} \Bigg ) \hspace{.2in} \text {\big ( by  (\ref {eq:gradK-refined}) \big )} \\
& = & C \left ( r_i^2 v_i(0)^{-\frac{2}{n-2}} + v_i(0) ^{2+2t_i-\frac{2}{n-2}} + 2 r_i \right ).
\end{eqnarray*}

\noindent
By (\ref {eq:limit2+2ti-2/n-2}) we know $\displaystyle \lim _{i \to \infty} \left (2+2t_i-\frac{2}{n-2}\right ) <0$, therefore $\displaystyle \lim _{i \to \infty}  v_i(0) ^{2+2t_i-\frac{2}{n-2}}=0$.  It follows from this and $\displaystyle \lim _{i \to \infty}  r_i^2 v_i(0)^{-\frac{2}{n-2}} = \displaystyle \lim _{i \to \infty} r_i=0$ that 
$$ \lim _{i \to \infty} v^2_i(0)\int _{B_{\sigma}} | X(K_i)| v_i^{p_i+1}\lambda _i^{\frac{2n}{n-2}}dy =0. $$
This completes the proof in case I.

\section{ Ruling Out Case II  }
\label{section:case2}

Now we consider Case II, which has been reduced to the following:
there is a sequence of functions $\{v_i\}$, each satisfies  
$$
\Delta _{g^{(i)}} v_i + K(\sigma _i y)v_i^{p_i}=0
$$
where $\sigma _i \to 0$ and $g^{(i)}(y)=g_{\alpha \beta}(\sigma _i y)dy^{\alpha}dy^{\beta}$. The sequence
$\{ v_i \}$ has isolated blow-up point(s) $\{0,...\}$. 

If $0$ is not a simple blow-up point, then we can do another rescaling and repeat the
argument in the previous section, with $r_i$ replaced by $r_i\sigma _i$, to get a contradiction.  Therefore $0$ must be a simple blow-up point for $\{v_i\}$.  Then we can still repeat the argument in the previous section, with $r_i$ replaced by $\sigma _i$. The only difference is in 
the expression of $h=\displaystyle \lim _{i \rightarrow \infty}
\frac{v_i(y)}{v_i(\bar{y})}$. 
 As shown in Section 7 of \cite{YY3}, because here $|y|^{\frac{2}{p_i-1}} \bar{v}_i(|y|)$ doesn't have a second critical point at $|y|=1$, we have a different expression of $h$: near $0$,
$$h(y) = c_1|y|^{2-n}+  A + O(|y|) $$ where $A$ is a positive constant.  This positive ``mass'' term $A >0 $ guarantees that the limit of the boundary term of the Pohozaev identity (\ref {eq:pohozaev-case1}) is negative, i.e., $$\lim _{i
\rightarrow \infty} \frac{1}{v^2_i(\bar{y})}\int _{\partial
B_{\sigma}} T_i(X, \nu _i) d \Sigma _i <0.$$  The other parts of the proof remain the same.  Therefore Case II can also be ruled out.

Thus we have finished the proof of Theorem \ref{thm:main}.

\section {Proof of Theorem \ref {thm:main-pos}}
\label {section:positive}

In this section we will prove Theorem \ref {thm:main-pos}.  There are many parallels between the proofs of Theorem \ref {thm:main} and Theorem \ref {thm:main-pos}.  Therefore we are going to emphasize the differences between the two proofs and omit the details of some of the steps if they can be obtained using essentially the same argument as in Theorem \ref {thm:main}.

By the standard elliptic theory, a bound on $\|u\|_{C^3(M)}$ can be easily obtained provided there is a uniform bound on $\|u\|_{C^0(M)}$.   Following from the Sobolev inequality and strong maximal principle, a uniform upper bound on $u$ would also imply a uniform lower bound away from $0$.  Therefore the main issue is to establish a uniform upper bound on all positive solutions $u$; again we prove this by contradiction.

Suppose this is not true, then there are sequences 
$\{ u_i \}$ and $ \{ p_i \} $ such that $$ \Delta _g u_i -c(n)R(g)u_i + Ku_i^{p_i}=0 \hspace{.3in} \text{ and } \hspace{.2in} \max _M u_i
\rightarrow \infty \,\, \text{ as } i \rightarrow \infty.$$  
By similar arguments as in the scalar-flat case, we can show that 
for fixed $\epsilon >0$ and $R >> 0
$  we can find $x_{1,i},..., x_{N(i),i}$ on $M$ for each function $u_i$ such that 
\begin{equation}
\text{ each } x_{j,i} \,\,
(1 \leq j \leq N(i) ) \text{ is a local maximum point of } u_i;
\end{equation}

\begin{equation}
\text{the balls } B_{\frac{R}{u_i(x_{j, i})^{\frac{p_i -1}{2}}}}(x_{j, i}) \text{ are disjoint};   
\end{equation}

\noindent
for coordinates $y=(y^1,...,y^n)$ such that $\frac{y}{u_i(x_{j,i})^{\frac{p_i-1}{2}}}$ are the conformally flat coordinates centered at $ x_{j,i}$,
\begin{equation}
\label{eq:iso-eq-pos}
\Bigg \| u_i(x_{j,i})^{-1}u_i \left (\frac
{y}{u_i(x_{j,i})^{\frac{p_i-1}{2}}}\right ) - \left (1+
\frac{K(x_{j, i})}{n(n-2)}|y|^2 \right )^{-\frac{n-2}{2}} \Bigg  \|_{C^2(B_{2R}(0))} <\epsilon;
\end{equation}
and
\begin{equation}
u_i(x) \leq C \left ( d_g(x,\{ x_{1,i}, ..., x_{N(i),i} \}) \right )^{-\frac{2}{p_i-1}} \hspace{.2in} \text{for a constant } C=C(\epsilon,R).       
\end{equation}
 
Let $\sigma _i= \min \{ d_g(x_{\alpha, i}, x_{\beta ,i}): \alpha \neq
\beta, 1 \leq \alpha, \beta \leq N(i) \}$. Without lost of generality we
can assume $\sigma _i= d_g(x_{1, i}, x_{2 ,i})$. As before there are two
possibilities.

\noindent
{\bf Case I}: \hspace{.05in}$\sigma _i \geq \varepsilon >0$.\\Then the points $x_{j,i}$ have isolated
limiting points $P_1, P_2, ...$, which are isolated blow-up points of $\{
u_i \}$.   

\noindent
{\bf Case II}: \hspace{.05in} $\sigma _i \rightarrow 0$.\\  Then we rescale the coordinates to make
the minimal distance 1: let $y=\sigma _i ^{-1}z$ where $z=(z^1,..., z^n)$ are the conformally flat coordinates centered at $x_{1,i}$. We also rescale the function by defining
$$v_i(y)=\sigma _i ^{\frac{2}{p_i-1}}u_i(\sigma _i y).$$
$v_i$ satisfies 
$$
\Delta _{g^{(i)}} v_i -c(n)R \left (g^{(i)}\right )v_i+ K_iv_i^{p_i}=0
$$
where $g^{(i)}(y)=g_{\alpha \beta}(\sigma _i y)dy^{\alpha}dy^{\beta}$, $R(g^{(i)})(y)=\sigma _i ^2 R(g)(\sigma _iy)$ and $K_i(y)=K(\sigma _iy)$.

We can prove as in Section 4 of \cite {YY3} that $0$ is an isolated blow-up point of $\{v_i\}$.

\subsection {Ruling Out Case I}
\label {subsection:case1-pos}

Now assume we are in Case I, i.e., all the blow-up points $\{P_1, P_2,...\}$ are isolated blow-up points.

\subsubsection{Simple Blow-up}
\label{subsubsection: simple-pos}

Next we need to study the behavior of the functions around simple blow-up points.  If any of the points, say $P_1$, is a simple blow-up point, then let $x_i$ be the local maximal point of $u_i$ such that $\displaystyle \lim _{i \to \infty} x_i=P_1$.  Let $z$ be the conformally flat coordinates centered at each $x_i$.  The next proposition is analogous to Proposition \ref {prop:simpleestimates}.

\begin{prop}
\label{prop:simpleestimates-pos}
There exist a constant $C$ independent of
$i$ and a radius $r_1 \leq \bar{r}$ (where $\bar{r}$ is defined as in Definition $\ref {defn:simple}$) such that
\begin{itemize}
\item if \,\, $0 \leq |z| \leq r_1$, \hspace{.1in} then $$u_i(z) \geq C u_i(x_i)\left (1+ \frac{K(x_i)}{n(n-2)}u_i(x_i)^{\frac{4}{n-2}}|z|^2 \right )^{-\frac{n-2}{2}}$$ 
\item  if \,\,  $0 \leq |z| \leq \frac{R}{u_i(x_i)^{\frac{p_i-1}{2}}}$, \hspace{.1in}
then $$u_i(z) \leq C u_i(x_i) \left (1+ \frac{K(x_i)}{n(n-2)}
u_i(x_i)^{p_i-1}|z|^2 \right )^{-\frac{n-2}{2}}$$ 
\item if \,\, $\frac{R}{u_i(x_i)^{\frac{p_i-1}{2}}} \leq |z| \leq r_1$,
\hspace{.1in} then \hspace{.1in} $u_i(z) \leq Cu_i(x_i)^{t _i}|z|^{-l_i}$\\
 where $l_i$, $t_i$ are so chosen that 
$\frac{(2n-1)(n-2)}{2n} < \displaystyle \lim_{i \to \infty} l_i < n-2$,
 and $t _i =1-\frac{(p_i-1)l_i}{2}$.
\end{itemize} 
\end{prop}

\noindent
Note that here $\displaystyle \lim_{i \to \infty} l_i$ is slightly different from that in Proposition \ref {prop:simpleestimates}; this modification is made to accommodate some adjustments (in a later part of the proof) that are related to $R(g)>0$.  However, the proof of this proposition is essentially the same as that of Proposition \ref {prop:simpleestimates}.  Therefore, in the proof below we will only point out the major steps and the few differences between the two proofs.  We refer the readers to the proof of Proposition \ref {prop:simpleestimates} for the details.

\pf
By (\ref {eq:iso-eq-pos}) when $0 \leq |z| \leq \frac{R}{u_i(x_i)^{\frac{p_i-1}{2}}}$,
\begin{eqnarray*}
& & (1-\epsilon) u_i(x_i)\left (1+ \frac{K(x_i)}{n(n-2)}
u_i(x_i)^{\frac{4}{n-2}}|z|^2\right )^{-\frac{n-2}{2}}\\
 & \leq &  u_i(z) \\
&  \leq & (1+\epsilon) u_i(x_i)\left (1+ \frac{K(x_i)}{n(n-2)}
u_i(x_i)^{p_i-1}|z|^2 \right )^{-\frac{n-2}{2}} .
\end{eqnarray*}

\noindent
So we only need to find the upper and lower bounds on $u_i(z)$ when $\frac{R}{u_i(x_i)^{\frac{p_i-1}{2}}} \leq |z| \leq r_1$.

\noindent
{\it The lower bound:}

\noindent
Let $G_i$ be the Green's function of $\Delta _g-c(n)R(g)$ which is singular at
$0$ and $G_i =0$ on $\partial B_{r_1}(x_i)$.  (Here the operator is different from the Laplacian operator which is used in the proof of Proposition \ref {prop:simpleestimates}).  By Lemma 9.2 in \cite {LZ}, 
there exist constants $C_1$ and $C_2$ independent of $i$ such that
 $$C_1 |z|^{2-n} \leq G_i(z) \leq C_2 |z|^{2-n}.$$ 

\noindent
There exists a constant $C$ independent of $i$, such that when $|z|=Ru_i(x_i)^{-\frac{p_i-1}{2}}$ and $|z|=r_1$, $Cu_i(x_i)^{-1} G_i(z) \leq u_i(z)$.

\noindent
Since
$$\Delta _g u_i - c(n)R(g)u_i \,\, = \,\, - Ku_i^{p_i+1} \,\, < \,\, 0$$
and 
$$\Delta _g G_i - c(n)R(g)G_i \,\, = \,\,  0,$$
we conclude by the maximal principle that $$u_i(z) \geq Cu_i(x_i)^{-1} G_i(z) \hspace {.2in} \text {when} \hspace {.1in} Ru_i(x_i)^{-\frac{p_i-1}{2}} \leq |z| \leq r_1. $$

\noindent
Finally because $G_i(z) \geq C_1 |z|^{2-n}$ and 
 $$u_i(x_i)^{-1}|z|^{2-n}
\geq C u_i(x_i) \left (1+ \frac{K(x_i)}{n(n-2)}
u_i(x_i)^{\frac{4}{n-2}}|z|^2 \right )^{-\frac{n-2}{2}}$$
for some constant $C$, we know
$$u_i(z) \geq C u_i(x_i) \left (1+ \frac{K(x_i)}{n(n-2)}
u_i(x_i)^{\frac{4}{n-2}}|z|^2 \right )^{-\frac{n-2}{2}}$$ 
when $\frac{R}{u_i(x_i)^{\frac{p_i-1}{2}}} \leq |z| \leq r_1$.

\noindent
{\it The upper bound:}

\noindent
Define $\mathcal{L}_i\varphi := \Delta _g \varphi -c(n)R(g) \varphi + K
u_i^{p_i-1}\varphi$. (The linear term is not in the $\mathcal{L}_i$ in the proof of Proposition \ref {prop:simpleestimates}).  By this definition $\mathcal{L}_iu_i=0$. Let $M_i=\displaystyle \max _{\partial B_{r_1}} u_i$ and $C_i=(1+\epsilon)
\left (\frac{K(x_i)}{n(n-2)} \right )^{-\frac{n-2}{2}}$. Note that $C_i$ is bounded above
and below by constants independent of $i$. Consider the function
$$M_i|z|^{-n+2+l_i}+C_iu_i(x_i)^{t _i}|z|^{-l_i}.$$

\noindent
On $\{|z|=r_1 \} \cup \{|z|=Ru_i(x_i)^{-\frac{p_i-1}{2}}\}$,$$u_i(z) \leq M_i|z|^{-n+2+l_i}+C_iu_i(x_i)^{t _i}|z|^{-l_i}.$$

\noindent
In the Euclidean metric, $\Delta |z|^{-l_i}=-l_i(n-2-l_i)|z|^{-l_i-2}
$ and $\Delta |z|^{-n+2+l_i}=-l_i(n-2-l_i)|z|^{-n+l_i}$.  Although here the metric $g$ may not be Euclidean, from the local coordinates expression of $\Delta_g$ it is easy to see that when $r_1$ is small enough, we can find a constant $C$ such that when $ |z| \leq r_1$,

\begin{equation*}
\Delta _g |z|^{-l_i} \leq -Cl_i(n-2-l_i)|z|^{-l_i-2}
\end{equation*}
and 
\begin{equation*}
\Delta _g |z|^{-n+2+l_i} \leq -Cl_i(n-2-l_i)|z|^{-n+l_i}. 
\end{equation*}
This implies
\begin{eqnarray*}
& & \mathcal{L}_i(C_iu_i(x_i)^{t _i}|z|^{-l_i}) \\
& = &  C_i u_i(x_i)^{t_i}\left ( \Delta _g|z|^{-l_i}-c(n)R(g)|z|^{-l_i}+Ku_i(z)^{p_i-1}|z|^{-l_i} \right)\\
& < & C_i u_i(x_i)^{t_i} \left (\Delta _g|z|^{-l_i}+Ku_i(z)^{p_i-1}|z|^{-l_i} \right ) \hspace {.3in} (\text {since} \hspace{.1in} R(g)>0)\\
& \leq & -Cl_i(n-2-l_i)u_i(x_i)^{t _i}|z|^{-l_i-2}+C'u_i(x_i)^{t _i}u_i(z)^{p_i-1}|z|^{-l_i}\\
& < & 0
\end{eqnarray*}
when $R$ is large enough, where the last inequality uses Lemma \ref {lemma:simpleharnack}, the simple blow-up property of $\{u_i\}$, and the fact that $l_i(n-2-l_i)$ is always
bounded below by some positive constant independent of $i$.

\noindent
Similarly, we can prove

$$\mathcal{L}_i(M_i|z|^{-n+2+l_i})<0. $$

\noindent
Therefore when $Ru_i(x_i)^{-\frac{p_i-1}{2}} \leq |z| \leq r_1$,
$$\mathcal{L}_i \left (M_i|z|^{-n+2+l_i}+C_iu_i(x_i)^{t
_i}|z|^{-l_i} \right ) <0, $$ and thus by the maximal principle  
$$
u_i(z)  \leq  M_i|z|^{-n+2+l_i}+C_iu_i(x_i)^{t_i}|z|^{-l_i}.
$$

\noindent
By Lemma \ref{lemma:simpleharnack} and the simple blow-up property of $\{ u_i \}$, for
$\frac{R}{u_i(x_i)^{\frac{p_i-1}{2}}} \leq   \theta \leq r_1$, 
\begin{eqnarray*}
M_i & \leq & C \theta^{\frac{2}{p_i-1}}\bar{u}_i(\theta) \\
& \leq & C\theta ^{\frac{2}{p_i-1}} \left (M_i\theta^{-n+2+l_i}+C_iu_i(x_i)^{t
_i} \theta ^{-l_i} \right )\\
 & = & C \theta ^{\frac{2}{p_i-1}-n+2+l_i}M_i + C \theta
^{\frac{2}{p_i-1}} \cdot C_iu_i(x_i)^{t_i} \theta ^{-l_i} 
\end{eqnarray*}
for some constant $C$ independent of $i$.

\noindent
Because $$ \displaystyle \lim_{i \to \infty}\left(\frac{2}{p_i-1}-n+2+l_i \right)\,\, = \,\, -\frac{n-2}{2}+  \displaystyle \lim_{i \rightarrow \infty} l_i \,\, > \,\,  -\frac{n-2}{2}+ \frac{(2n-1)(n-2)}{2n} \,\, > \,\, 0$$ 

\noindent
and $\frac{R}{u_i(x_i)^{\frac{p_i-1}{2}}} \to 0$, we can choose $\theta$ small enough (fixed and independent of $i$) to
absorb the first term on the right hand side of the above inequality
into the left hand side to get $M_i \leq  C u_i(x_i)^{t _i}$.  Therefore
\begin{eqnarray*}
u_i(z) & \leq &  M_i|z|^{-n+2+l_i}+C_iu_i(x_i)^{t_i}|z|^{-l_i}\\
& \leq & M_i|z|^{-l_i}+C_iu_i(x_i)^{t_i}|z|^{-l_i}\\ 
& \leq &  C u_i(x_i)^{t _i}|z|^{-l_i}
\end{eqnarray*} 
     
\stop

The following technical lemma is parallel to Lemma \ref {lemma:intv^p+1y^k}.  Note that because of the modification of $\displaystyle \lim _{i \to \infty} l_i$ we are able to have the estimate up to $\kappa = n-1$.

\begin{lemma}
\label {lemma:intu^p+1z^k-pos}
When $\sigma <r_1$ and $0 \leq \kappa \leq n-1$, $$\int _{|z|\leq \sigma} |z|^{\kappa}u_i(z)^{p_i+1}dz \leq C u_i(x_i)^{-\frac{2\kappa}{n-2}+\frac{n-2+\kappa}{2}\delta_i},$$ where $C$ is independent of $i$ and $r_1$ is defined as in Proposition $\ref {prop:simpleestimates-pos}$.
\end{lemma}

\pf
\noindent
By Proposition \ref{prop:simpleestimates-pos}
\begin{eqnarray*}
 \int _{|z| \leq \frac{R}{u_i(x_i)^{\frac{p_i-1}{2}}} } |z|^{\kappa}
u_i(z)^{p_i+1} dz & \leq & C
u_i(x_i)^{p_i+1}  \int _{|z| \leq \frac{R}{u_i(x_i)^{\frac{p_i-1}{2}}}
} |z|^{\kappa} dz   \\
& \leq & C u_i(x_i)^{p_i+1-\frac{(n+\kappa)(p_i-1)}{2}}\\
& = & C u_i(x_i) ^{-\frac{2\kappa}{n-2} + \frac{n-2+\kappa}{2}\delta _i}.
\end{eqnarray*}

\noindent
By our choice of $\l_i$ 
\begin{eqnarray*}
\lim _{i \rightarrow \infty} \big ( n+\kappa-l_i(p_i+1) \big ) &  = &  n+\kappa-\frac{2n}{n-2}\lim_{i \to \infty} l_i \\
& < & n+ \kappa-\frac{2n}{n-2}\cdot \frac{(2n-1)(n-2)}{2n}\\
& \leq & n+(n-1)-(2n-1)\\
& = & 0.
\end{eqnarray*}
Therefore
\begin{eqnarray*}
\int _{\frac{R}{u_i(x_i)^{\frac{p_i-1}{2}}} \leq |z| \leq \sigma } |z|^{\kappa}u_i(z)^{p_i+1} dz & \leq & C \int
_{\frac{R}{u_i(x_i)^{\frac{p_i-1}{2}}} \leq |z| \leq \sigma } |z|^{\kappa}
\left ( u_i(x_i)^{t _i}|z|^{-l_i} \right )^{p_i+1} dz \\
            & \leq & C u_i(x_i)^{ t _i(p_i+1) -\frac{p_i-1}{2}\left
(n + \kappa -l_i(p_i+1) \right )} \\
& = & C u_i(x_i)^{p_i+1-\frac{(n+\kappa)(p_i-1)}{2}} \hspace{.2in}(\text{by the definition of }  t _i)\\
 & = & C u_i(x_i) ^{-\frac{2\kappa}{n-2} +\frac{n-2+\kappa}{2}\delta_i}.
\end{eqnarray*}
Thus
\begin{equation*}
\int _{|z| \leq \sigma } |z|^{\kappa} u_i(z)^{p_i+1} dz \leq C u_i(x_i)
^{-\frac{2\kappa}{n-2} + \frac{n-2+\kappa}{2}\delta _i}.
\end{equation*}

\stop

Let $\delta _i := \frac{n+2}{n-2}-p_i$.  Since the background metric $g$ is locally conformally flat, we can write it locally as $\lambda (z) ^{\frac{4}{n-2}} dz^2.$  Let $\sigma < r_1 $.  As in the scalar-flat case, we need to use the Pohozaev identity: for a conformal Killing field $X$ on $B_{\sigma}(x_i)$,  
\begin{equation}
\label{eq:pohozaev-positive} 
\frac{n-2}{2n} \int _{B_{\sigma}} X(R_i) dv_{g_i} = \int _{\partial
B_{\sigma}} T_i(X, \nu _i) d \Sigma _i
\end{equation}

\noindent
where 
\begin{eqnarray*}
g_i &  = &  u_i^{\frac{4}{n-2}}g \,\,  = \,\, (\lambda u_i)^{\frac{4}{n-2}}dz^2, \\
R_i &  = & R(g_i) \,\, = \,\, c(n)^{-1} Ku_i^{- \delta _i},\\
 dv_{g_i} & = &  (\lambda u_i)^{\frac{2n}{n-2}}dz  ,\\
\nu _i & = & (\lambda u_i)^{-\frac{2}{n-2}}\sigma^{-1}\sum _{j} z^j
\frac{\partial}{\partial z^j} \\
& &  \text { is the unit outer
normal vector on } \partial B_{\sigma} \text{ with respect to } g_i,\\
d\Sigma _i & = & (\lambda u_i)^{\frac{2(n-1)}{n-2}} d \Sigma _{\sigma}
\\
& & 
\text{ where } d \Sigma _{\sigma} \text { is the surface element of the
standard } S^{n-1}(\sigma),\\
T_i & = & \Ric(g_i)-n^{-1}R(g_i)g_i \hspace{.2in} \text{ is the traceless Ricci tensor with respect to } g_i.
\end{eqnarray*}

\noindent
$T_i$ can also be expressed as  $$(n-2) (\lambda u_i)^{\frac{2}{n-2}} \left(\Hess \left ( (\lambda u_i)^{-\frac{2}{n-2}} \right )-\frac{1}{n}\Delta \left ((\lambda u_i)^{-\frac{2}{n-2}} \right ) dz^2 \right )$$
\noindent
where $\Hess$ and $\Delta$ are taken with respect to the Euclidean metric $dz^2$.

Now we choose $X=\displaystyle \sum _{j=1}^{n} z^j \frac {\partial}{\partial z^j}$.  By an argument which is almost identical to that in the proof of Proposition \ref {prop:v^delta} we know

\begin{eqnarray*}
& & \frac{n}{2(n-1)}\left(1+ \frac{\delta _i}{p_i+1}\right) \int  _{B_{\sigma}}  |z| \frac {\partial K}{\partial r} u_i^{p_i+1}\lambda ^{\frac{2n}{n-2}}dz  + \frac{n}{2(n-1)}\frac{\delta _i}{p_i+1}n  \int
_{ B_{\sigma}} Ku_i^{p_i+1}\lambda ^{\frac{2n}{n-2}} dz \\
& \leq  & Cu_i(x_i)^{2t_i}+C\delta _i u_i(x_i)^{t_i(p_i+1)} + C\delta_i u_i(x_i)^{-\frac{2}{n-2}+\frac{n-1}{2}\delta_i}. 
\end{eqnarray*}

\noindent
Since
\begin{eqnarray}
\label{eq:lower-int-u^p+1-pos}
 \int _{ B_{\sigma}} u_i^{p_i+1} dz & > & \int _{|z| \leq
\frac{R}{u_i(x_i)^{\frac{p_i-1}{2}} } } u_i^{p_i+1} dz \nonumber \\
& \geq & C u_i(x_i)^{p_i+1 - \frac{n}{2}(p_i-1)} \hspace{.3in} \text {(by Proposition } \ref {prop:simpleestimates-pos})   \nonumber \\
& = & C u_i(x_i)^{ \frac{n-2}{2} \delta _i} \nonumber \\
& \geq & C,
\end{eqnarray}
as before we can argue that 
$$\delta _i \leq  C\bigg | \int  _{B_{\sigma}}   \frac {\partial K}{\partial r}|z| u_i^{p_i+1}\lambda^{\frac{2n}{n-2}}dz \bigg | +  C u_i(x_i)^{2t _i}  + C \delta_i u_i(x_i)^{t_i(p_i+1)}   + C\delta_i u_i(x_i)^{-\frac{2}{n-2}+\frac{n-1}{2}\delta_i} .$$

\noindent
Then because
\begin{equation}
\label {eq: lim-ti-pos}
\lim _{i \to \infty} t_i \,\, = \,\, 1 - \frac {2}{n-2} \lim  _{i \to \infty} l_i \,\, < \,\, 1- \frac {2}{n-2} \cdot \frac{(2n-1)(n-2)}{2n} \,\, = \,\, \frac {1-n}{n} \,\, < \,\, 0,
\end{equation}
the last two terms on the right hand side can be absorbed into the left hand side, so we have

\begin{equation}
\label {eq:delta-coarse-pos} 
\delta _i \leq C \left (\bigg | \int _{B_{\sigma}} \frac{\partial K}{\partial r} |z|u_i^{p_i+1} \lambda (z)^{\frac{2n}{n-2}} dz \bigg |+u_i(x_i)^{2t_i}  \right ).
\end{equation}

\noindent
By Lemma \ref {lemma:intu^p+1z^k-pos} this implies 
\begin{equation}
\label {eq:u^delta-pos}
\lim _{i \to \infty} u_i(x_i)^{\delta _i} =1,
\end{equation}
which is parallel to Proposition \ref {prop:v^delta}; and we also have a preliminary estimate for $\delta _i$:
\begin{equation}
\label {eq:delta-prelim-pos}
\delta _i \leq C \left ( u_i(x_i)^{2t_i} + u_i(x_i)^{-\frac{2}{n-2}} \right ).
\end{equation}

{\it Now suppose the blow-up points $\{P_1, P_2, ...\}$ are all simple blow-up points.}  Choose a point $P \in \partial B_{\frac{r_1}{2}}(P_1)$, by Proposition \ref {prop:simpleestimates-pos} we know $u_i(P) \to 0$ as $i \to \infty$.  Let $\Omega$ be any compact subset of $M \setminus \{P_1, P_2, ...\}$ containing $P$.  By Definition \ref {defn:iso}, $u_i$ is bounded above on $\Omega$ by some constant $C$ independent of $i$ (although it may depend on $\Omega$), thus on $\Omega$ we have the standard Harnack inequality.  Therefore 
$$ \max _{\Omega} \frac{u_i}{u_i(P)} \,\, \leq \,\, C \min  _{\Omega} \frac{u_i}{u_i(P)} \,\, \leq \,\, C \frac{u_i(P)}{u_i(P)} \,\, = \,\, C.  $$
Since $u_i$ satisfies (\ref {eq:yamabe}),
$$ \Delta _g \left ( \frac{u_i}{u_i(P)}  \right ) -c(n)R(g)\frac{u_i}{u_i(P)}+u_i(P)^{p_i-1}K \left ( \frac{u_i}{u_i(P)}  \right )^{p_i}=0.  $$
Then by the standard elliptic theory, $ \frac{u_i}{u_i(P)}$ converges in $C^2$-norm on $\Omega$ to some function $G \geq 0$ which satisfies $\Delta _g G-c(n)R(g)G=0$ on $\Omega$.  Because $\Omega$ is arbitrary, $G$ satisfies $\Delta _g G-c(n)R(g)G=0$ on $M \setminus \{P_1, P_2, ...\}$.  Since $R(g)>0$, $G$ must be singular at one or more of the points $\{P_1, P_2, ...\}$.  Suppose it is singular at $P_1,...,P_k$, it follows that $G$ is a linear combination of the positive fundamental solutions $G_{\gamma}$ with poles at $P_{\gamma}$ for $\gamma = 1,...,k$, i.e., there exist positive constants $a_1,..., a_k$ such that $G=\displaystyle \sum _{\gamma =1}^{k} a_{\gamma}G_{\gamma}$.

{\it This is precisely the key difference between the scalar-flat and the scalar-positive cases.  Recall that when $R(g) \equiv 0$, we used a removable singularity theorem for harmonic functions to prove that the isolated blow-up points cannot all be simple (Section 6 of \cite{YY3}).  Here because $R(g)>0$, we will need to do more work to show that. }

Next we apply the Pohozaev identity (\ref {eq:pohozaev-positive}) to $X=\frac{\partial}{\partial z^1}$.  As in the scalar-flat case, direct computation shows that the boundary term is equal to 

\begin{eqnarray*}
& & (n-2)
\int _{\partial B_{\sigma}} 
 \sum _{j} \frac{z^j}{\sigma}    \left( -\frac {2 }{n-2} (\lambda u_i) 
\frac{\partial ^2 (\lambda u_i)  }{\partial z^1 \partial z^j}  
 +  \frac{2n}{(n-2)^2}
\frac{\partial (\lambda u_i)}{\partial z^1} \frac{\partial( \lambda u_i)}{\partial z^j}  \right) \\ 
& & - \frac{z^1}{\sigma} \sum _{j}\left( -\frac{2 }{n(n-2)} (\lambda u_i)  \frac{\partial ^2  (\lambda u_i) }{(\partial
z^j)^2}  + \frac{2}{(n-2)^2} \left(\frac{\partial (\lambda u_i)}{\partial z^j} \right )^2 \right )
d \Sigma _{\sigma}, 
\end{eqnarray*}
and it decays in the rate of $u_i(x_i)^{2 t _i}$.

\noindent
The interior term 

\begin{eqnarray*}
 \frac{n-2}{2n} \int _{B_{\sigma}} \frac{\partial}{\partial
z^1}(R_i) dv_{g_i} & = & \frac{n-2}{2n} c(n)^{-1} \int _{B_{\sigma}} \left (1+
\frac{\delta_i}{p_i+1} \right )\lambda^{\frac{2n}{n-2}}u_i^{p_i+1}\frac{\partial K}{\partial z^1} dz \\
& & +\frac{n-2}{2n}c(n)^{-1} \int _{B_{\sigma}} \frac{\delta_i}{p_i+1}
K u_i^{p_i+1} \frac {\partial \lambda ^{\frac{2n}{n-2}}}{\partial z^1} dz\\
& & - \frac{n-2}{2n} c(n)^{-1}\frac{\delta_i}{p_i+1}\int _{\partial B_{\sigma}}
\lambda^{\frac{2n}{n-2}}Ku_i^{p_i+1} \frac{z^1}{\sigma}d \Sigma _{\sigma}.
\end{eqnarray*}

\noindent
By Proposition \ref {prop:simpleestimates-pos}, Lemma \ref {lemma:intu^p+1z^k-pos} and (\ref {eq:u^delta-pos}), the second term is bounded by 
 $$C\delta_i u_i(x_i)^{\frac{n-2}{2}\delta_i}\leq C\delta _i,$$
and the last term is
bounded by $$C\delta _i u_i(x_i)^{t_i(p_i+1)} \leq C\delta _iu_i(x_i)^{2t _i}
.$$

\noindent
Thus we have a bound on the first term: 
\begin{eqnarray*}
\frac{n-2}{2n} c(n)^{-1} \int _{B_{\sigma}} \left(1+
\frac{\delta_i}{p_i+1}\right)\lambda^
{\frac{2n}{n-2}}u_i^{p_i+1}\frac{\partial K}{\partial z^1} dz & \leq &
C(u_i(x_i)^{2 t _i}+\delta _i u_i(x_i)^{2t_i}+ \delta _i )\\
& \leq & C \left ( u_i(x_i)^{2 t _i} + \delta _i\right ).
\end{eqnarray*}

\noindent
This shows that 
\begin{equation}
\label{eq:gradK-coarse-pos}
\bigg| \int _{B_{\sigma}} \lambda ^{\frac{2n}{n-2}}u_i^{p_i+1}\frac{\partial K}{\partial z^1} dz \bigg| \leq C \left (  u_i(x_i)^{2 t _i} + \delta _i \right ).
\end{equation}

\noindent
By the Taylor expansion,
$$\frac{\partial K}{\partial z^1}(z)= \frac{\partial K}{\partial
z^1}(0) + \nabla \left(\frac{\partial K}{\partial
z^1} \right )(\varsigma) \cdot z
\hspace{.3in} \text{ for some } |\varsigma| \leq |z|.$$

\noindent
By Lemma \ref{lemma:intu^p+1z^k-pos} and (\ref {eq:u^delta-pos}),

\begin{eqnarray*}
 \int _{B_{\sigma}} \lambda ^{\frac{2n}{n-2}}u_i^{p_i+1} \Bigg | \nabla \left (\frac{\partial K}{\partial z^1} \right )(\varsigma) \cdot z \Bigg | dz & \leq &  C \int _{B_{\sigma}} u_i^{p_i+1}|z|dz\\
& \leq & Cu_i(x_i)^{-\frac{2}{n-2}}.  
\end{eqnarray*}
 
\noindent
Together with (\ref {eq:lower-int-u^p+1-pos}) and (\ref {eq:gradK-coarse-pos}), this shows that
\begin{eqnarray*}
\bigg | \frac{\partial K}{\partial z^1}(x_i) \bigg | \,\, = \,\, 
\bigg | \frac{\partial K}{\partial z^1}(0) \bigg | & \leq & C \left (u_i(x_i)^{-\frac{2}{n-2}} + u_i(x_i)^{2 t _i}+\delta _i \right )\\
& \leq & C \left (u_i(x_i)^{-\frac{2}{n-2}} + u_i(x_i)^{2 t _i} \right ) \hspace{.3in} \text {by } (\ref {eq:delta-prelim-pos}).
\end{eqnarray*}

\noindent
The same estimate holds for $\big |\frac{\partial K}{\partial z^j}(x_i) \big |, \,\, j=2,...,n $ as well, so we know $|\nabla K (P_1)| = \displaystyle \lim _{i \to \infty} | \nabla K (x_i)| =0$.  That is, the blow-up point $P_1$ is a critical point of $K$.

In the next step we once again study the Pohozaev identity with $X=\displaystyle \sum _{j}z^j\frac {\partial}{\partial z^j}$.  We divide both sides of it by $u_i^2(P)$, so it becomes

\begin{equation}
\label{eq:pohozaev-divide-pos}
\frac{n-2}{2n} \frac{1}{u_i^2(P)}\int _{B_{\sigma}(x_i)} X(R_i) dv _{g_i} =\frac{1}{u_i^2(P)} \int _{\partial
B_{\sigma}(x_i)} T_i(X, \nu _i) d \Sigma _i.
\end{equation}

The right hand side (boundary term) is 

\begin{eqnarray*}
 & & \frac{1}{u_i^2(P)}\int _{\partial B_{\sigma}(x_i)} T_i(X, \nu _i) d \Sigma _i\\ & = & \frac{1}{u_i^2(P)}\int _{\partial B_{\sigma}(x_i)} 
\bigg [ \Ric 
  \left ( 
    \left (\lambda u_i \right )^{\frac{4}{n-2}}dz \otimes dz \right )  \\
& & - n^{-1} R 
  \left( \left(\lambda u_i \right)^{\frac{4}{n-2}}dz \otimes dz \right) 
  \left( \lambda u_i \right )^{\frac{4}{n-2}}dz \otimes dz
\bigg ](X, \nu_0) (\lambda u_i)^2 d \Sigma _{\sigma}  \\
& = &  \int _{\partial B_{\sigma}(x_i)} 
\left( \frac{\lambda u_i}{u_i(P)}\right)^2 
\Bigg [ \Ric
  \left ( 
     \left(\frac{\lambda u_i}{u_i(P)}\right)^{\frac{4}{n-2}} dz \otimes dz
  \right) \\
  & & - n^{-1} R 
  \left ( 
    \left(\frac{\lambda u_i}{u_i(P)}\right)^{\frac{4}{n-2}}dz \otimes dz
  \right) 
  \left( \frac{\lambda u_i}{u_i(P)} \right ) 
  ^{\frac{4}{n-2}}dz \otimes dz
\Bigg ] (X, \nu _0)    d \Sigma_{\sigma} 
\end{eqnarray*} 
 where $\nu _0 =\sigma ^{-1} \displaystyle \sum _{j}z^j \frac{\partial}{\partial z^j}$ is the unit outer normal on $\partial B_{\sigma} (x_i)$ with respect to 
the Euclidean metric $dz \otimes dz$. 

\noindent
Recall that on $B_{\sigma} (P_1) \setminus \{ P_1\}$, $\frac{u_i}{u_i(P)} \to G$ as $i \to \infty$, so the boundary term converges to

$$ \int _{\partial
B_{\sigma}(P_1)} (\lambda G)^2  \bigg (\Ric \left ((\lambda G)^{\frac{4}{n-2}}dz ^2 \right )  - n^{-1} R \left ((\lambda G)^{\frac{4}{n-2}} dz ^2 \right )(\lambda G)^{\frac{4}{n-2}}  dz^2 \bigg )(X, \nu _0)    d \Sigma _{\sigma}, $$ which can be expressed as


\begin{equation}
\label{eq:bdry-term-pos}
(n-2)\sigma^{-1} \int _{\partial
B_{\sigma}(P_1)} (\lambda G)^{\frac{2(n-1)}{n-2}} \cdot \left [ \Hess \left
((\lambda G)^{-\frac{2}{n-2}} \right )(X,X)-\frac{1}{n}\Delta \left ((\lambda G)^{-\frac{2}{n-2}} \right )\sigma ^2  \right ] \,\,  d \Sigma _{\sigma}. 
\end{equation}

\noindent
Since $\Delta _g G-c(n)R(g)G=0$ on $B_{\sigma}(P_1) \setminus \{ P_1 \}$, we know $G^{\frac{4}{n-2}}g=(\lambda G)^{\frac{4}{n-2}}dz^2$ has zero scalar curvature.  This implies that $\lambda G$ is a positive Euclidean harmonic function on $B_{\sigma} (0) \setminus \{ 0 \}$ which is singular at $0$.  Therefore  $\lambda G$ has the expression
$$ (\lambda G)(z)=a_1|z|^{2-n}+A+h(z) $$
where $h(z)$ is a harmonic function with $h(0)=0$.  Furthermore, the fundamental solution $G_1$ satisfies $$(\lambda G_1)(z)=|z|^{2-n}+E(P_1)+O(|z|).$$  Here $E(P_1)$ is the energy at $P_1$, and by the Positive Mass Theorem \cite {SY3} $E(P_1)>0$ since $(M,g)$ is not conformally equivalent to $S^n$.  Then because $G \geq a_1G_1$, we know that $A \geq a_1E(P_1) >0$.

\noindent
Next we calculate (\ref {eq:bdry-term-pos}).

\begin{eqnarray*}
(\lambda G)^{-\frac{2}{n-2}} & = & \left( a_1|z|^{2-n}+A+O(|z|) \right )^{-\frac{2}{n-2}}\\
& = &  a_1^{-\frac{2}{n-2}}|z|^2 - \frac{2}{n-2}Aa_1^{-\frac{n}{n-2}}|z|^{n} + O \left ( |z|^{2n-2} \right ).
\end{eqnarray*}

\noindent
Since
$$\Hess \left( a_1^{-\frac{2}{n-2}}|z|^2 - \frac{2}{n-2}Aa_1^{-\frac{n}{n-2}}|z|^{n} \right )(X,X) = 2a_1^{-\frac{2}{n-2}}|z|^2- \frac{2(n^2-n)}{n-2}Aa_1^{-\frac{n}{n-2}}|z|^{n},$$

\noindent
we have

$$\Hess \left ( (\lambda G)^{-\frac {2}{n-2} } \right ) (X, X)= 2a_1^{-\frac{2}{n-2}}|z|^2- \frac{2(n^2-n)}{n-2}Aa_1^{-\frac{n}{n-2}}|z|^{n}+O(|z|^{2n-2});$$

\noindent
and because
$$\Delta \left (  a_1^{-\frac{2}{n-2}}|z|^2 - \frac{2}{n-2}Aa_1^{-\frac{n}{n-2}}|z|^{n} \right )=2n a_1^{-\frac{2}{n-2}}-\frac{2(2n^2-2n)}{n-2}Aa_1^{-\frac{n}{n-2}}|z|^{n-2}, $$

\noindent
we have

$$\frac{1}{n}\Delta \left ((\lambda G)^{-\frac{2}{n-2}} \right )\sigma ^2 = 2 a_1^{-\frac{2}{n-2}}\sigma ^2-\frac{2(2n-2)}{n-2}Aa_1^{-\frac{n}{n-2}}|z|^{n-2}\sigma ^2 + O(|z|^{2n-4})\sigma ^2 .$$

\noindent
Therefore on $\partial B_{\sigma}(P_1)$,
$$ \Hess \left ((\lambda G)^{-\frac{2}{n-2}} \right )(X,X)
-\frac{1}{n}\Delta \left ((\lambda G)^{-\frac{2}{n-2}} \right )\sigma ^2 =  -2(n-1)Aa_1^{-\frac{n}{n-2}} \sigma ^n + O \left (\sigma ^{2n-2} \right )  . $$    
We also know 

\begin{eqnarray*}
(\lambda G)^{\frac{2(n-1)}{n-2}} & = & \left( a_1|z|^{2-n}+A+O(|z|) \right )^{\frac{2(n-1)}{n-2}}\\
& = &  a_1^{\frac{2(n-1)}{n-2}}|z|^{-2(n-1)} \left (1+  \frac{2(n-1)}{n-2}\frac{A}{a_1}|z|^{n-2} + O ( |z|^{2n-2}) \right ).
\end{eqnarray*}

\noindent
Thus (\ref{eq:bdry-term-pos}) is equal to 
\begin{eqnarray}
\label{eq:bdry-term-sign-pos}
 &  & (n-2) \sigma^{-1} \cdot a_1^{\frac{2(n-1)}{n-2}}\sigma^{-2(n-1)} \bigg (1+  \frac{2(n-1)}{n-2}\frac{A}{a_1}\sigma^{n-2} \nonumber \\
& & + O ( \sigma^{2n-2}) \bigg )  \cdot \left (-2(n-1)Aa_1^{-\frac{n}{n-2}} \sigma ^n + O (\sigma ^{2n-2}) \right ) \sigma^{n-1}\nonumber \\
& = &  -2(n-1)(n-2)Aa_1 + O(\sigma ^{n-2}) \nonumber\\
& < & 0 \nonumber
\end{eqnarray}
when $\sigma$ is sufficiently small, since $A>0$ and $a_1>0$. 

On the other hand, the left hand side (interior term) of (\ref{eq:pohozaev-divide-pos}) is 
$$
\frac{n-2}{2n}c(n)^{-1} \frac{1}{u^2_i(P)}\int _{B_{\sigma}(x_i)}
X(Ku_i^{-\delta _i})(\lambda u_i)^{\frac{2n}{n-2}}dz.
$$
Using the divergence theorem we can write
\begin{eqnarray*}
& & \frac{1}{u^2_i(P)}\int _{B_{\sigma}(x_i)}X(Ku_i^{-\delta _i })(\lambda u_i)^{\frac{2n}{n-2}}dz  \\
 & = & \frac{1}{u^2_i(P)}\int _{B_{\sigma}(x_i)}
X(K)u_i^{p_i +1 }\lambda ^{\frac{2n}{n-2}}dz  -\frac{\delta _i}{p_i+1}\frac{\sigma}{u^2_i(P)} \int
_{\partial B_{\sigma}(x_i)} K u_i^{p_i+1}\lambda ^{\frac{2n}{n-2}} d \Sigma _{\sigma} \\
& &  + \frac{\delta _i}{p_i+1}\frac{1}{u^2_i(P)} \int _{B_{\sigma}(x_i)} Ku_i^{p_i+1} \lambda ^{\frac{2n}{n-2}} \left (n + X(\ln K)+ \frac{2n}{n-2}X(\ln \lambda ) \right ) dz.
\end{eqnarray*}

\noindent
The second term
\begin{eqnarray*}
 -\frac{\delta _i\sigma}{p_i+1}\frac{1}{u^2_i(P)} \int
_{\partial B_{\sigma}(x_i)} K u_i^{p_i+1}\lambda ^{\frac{2n}{n-2}} d \Sigma _{\sigma} &  = &  -\frac{\delta _i\sigma}{p_i+1} \int
_{\partial B_{\sigma}(x_i)} K  \left ( \frac {u_i}{u_i(P)} \right )^2 u_i^{p_i-1}  \lambda ^{\frac{2n}{n-2}} d \Sigma_{\sigma} \\
& \rightarrow & 0,
\end{eqnarray*}
since $\frac{u_i}{u_i(P)} \rightarrow G$ and $u_i \rightarrow 0$ uniformly on $B_{2 \sigma}(P_1)\setminus B_{\frac{\sigma}{2}}(P_1)$.

\noindent
Since $X=r\frac{\partial}{\partial r}$ and $\frac{\partial}{\partial
r}(\ln K)$, $\frac{\partial}{\partial r}(\ln \lambda )$ are uniformly bounded, we can choose $\sigma$ small (independent of $i$) to make $n + X(\ln K )+ \frac{2n}{n-2}X(\ln \lambda )>0$.  Thus the limit of the last term is greater than or equal to $0$.  

\noindent
We claim that the first term $\frac{1}{u^2_i(P)}\int _{B_{\sigma}(x_i)}
X(K)u_i^{p_i +1 }\lambda ^{\frac{2n}{n-2}}dz \to 0$.  It follows from Proposition \ref {prop:simpleestimates-pos} that $u_i(P)\geq Cu_i(x_i)^{-1}$, thus to prove this limit it suffices to show that  
\begin{equation}
\label{eq:key-X(K)-pos}
\lim _{i \rightarrow \infty} u^2_i(x_i)\int _{B_{\sigma}(x_i)}
X(K)u_i^{p_i+1}\lambda ^{\frac{2n}{n-2}}dz =0.
\end{equation}

We write $ X(K)= r \frac {\partial K}{\partial r}=\displaystyle \sum_{j=1}^nz^j\frac {\partial K}{\partial z^j}$.  Since the coordinates are centered at $x_i$, for each $j=1,...,n$,
\begin{eqnarray*}
\frac {\partial K}{\partial z^j}(z) & = & \frac {\partial K}{\partial z^j}(x_i) + \sum _{|\beta|=1} \frac{\partial ^{\beta}}{\partial z^{\beta}}\frac {\partial K}{\partial z^j}(x_i)z^{\beta}+\frac{1}{2!}\sum _{|\beta|=2} \frac{\partial ^{\beta}}{\partial z^{\beta}}\frac {\partial K}{\partial z^j}(x_i)z^{\beta}+\cdot \cdot \cdot \\
& & + \frac{1}{(n-3)!}\sum _{|\beta|=n-3} \frac{\partial ^{\beta}}{\partial z^{\beta}}\frac {\partial K}{\partial z^j}(x_i)z^{\beta}+\frac{1}{(n-2)!}\sum _{|\beta|=n-2} \frac{\partial ^{\beta}}{\partial z^{\beta}}\frac {\partial K}{\partial z^j}(\varsigma)z^{\beta}
\end{eqnarray*}
where $|\varsigma|\leq |z|$. Therefore

\begin{eqnarray*}
& &   \int  _{B_{\sigma}} \bigg | r \frac {\partial K}{\partial r} \bigg | u_i^{p_i+1}\lambda ^{\frac{2n}{n-2}}dz  \\
 & \leq & C \Bigg ( \int _{B_{\sigma}} \bigg|  \frac {\partial K}{\partial z^j}(x_i)\bigg| |z|u_i^{p_i+1}dz + \sum _{|\beta|=1}^{n-3} \int _{B_{\sigma}} \bigg | \frac{\partial ^{\beta}}{\partial z^{\beta}}\frac {\partial K}{\partial z^j}(x_i)\bigg| |z|^{|\beta|+1} u_i^{p_i+1}dz\\
& & + \sum _{|\beta|=n-2} \int _{B_{\sigma}} \bigg | \frac{\partial ^{\beta}}{\partial z^{\beta}}\frac {\partial K}{\partial z^j}(\varsigma)\bigg| |z|^{n-1} u_i^{p_i+1}dz \Bigg ).
\end{eqnarray*}

\noindent
By Lemma \ref{lemma:intu^p+1z^k-pos} and (\ref {eq:u^delta-pos}), $$\int _{B_{\sigma}} \bigg |  \frac {\partial K}{\partial z^j}(x_i)\bigg | |z|u_i^{p_i+1}dz \leq C   |\nabla K(x_i)  | u_i(x_i)^{-\frac{2}{n-2}}$$ 
and 
$$\sum _{|\beta|=n-2} \int _{B_{\sigma}} \bigg | \frac{\partial ^{\beta}}{\partial z^{\beta}}\frac {\partial K}{\partial z^j}(\varsigma)\bigg| |z|^{n-1} u_i^{p_i+1}dz \leq Cu_i(x_i)^{-\frac{2(n-1)}{n-2}}. $$

\noindent
In addition, because $K$ satisfies the flatness condition $(**)$, as in the scalar-flat case we can show that
$$\bigg | \frac{\partial ^{\alpha} K}{\partial z^{\alpha}}(x_i) \bigg | \leq C | \nabla K(x_i) |^{\frac{n-1-|\alpha|}{n-2}} $$ when $  2\leq |\alpha|\leq n-2 .$
\noindent
Thus for any $1 \leq |\beta| \leq n-3$, 

\begin{eqnarray}
& & \int _{B_{\sigma}} \bigg | \frac{\partial ^{\beta}}{\partial z^{\beta}}\frac {\partial K}{\partial z^j}(x_i)\bigg| |z|^{|\beta|+1} u_i^{p_i+1}dz \nonumber \\
 & \leq & C \int _{B_{\sigma}}  | \nabla K(x_i) |^{\frac{n-1-(|\beta|+1)}{n-2}}   |z|^{|\beta|+1} u_i^{p_i+1}dz \nonumber \\
& = & C \int _{B_{\sigma}} | \nabla K(x_i) |^{\frac{n-2-|\beta|}{n-2}}|z|^{|\beta|}\cdot|z|u_i^{p_i+1} dz \nonumber \\
& \leq & C \int _{B_{\sigma}}\left ( | \nabla K(x_i) |^{\frac{n-2-|\beta|}{n-2}\cdot \frac{n-2}{n-2-|\beta|}}+ |z|^{|\beta|\cdot \frac{n-2}{|\beta|} } \right )\cdot|z|u_i^{p_i+1} dz \nonumber \\
& & \text{ (by Young's Inequality)} \nonumber \\
& = & C \left( \int _{B_{\sigma}}  | \nabla K(x_i) | \cdot |z|u_i^{p_i+1} dz + \int _{B_{\sigma}} |z|^{n-1}u_i^{p_i+1} dz  \right) \nonumber \\
& \leq &  C | \nabla K(x_i) | u_i(x_i)^{-\frac{2}{n-2}}+ C u_i(x_i)^{-\frac{2(n-1)}{n-2}}. \nonumber
\end{eqnarray}
Thus 
\begin{equation}
\label {eq:int-rdK/dr-u^p+1-pos}
 \int  _{B_{\sigma}} \bigg | r \frac {\partial K}{\partial r} \bigg | u_i^{p_i+1} \lambda ^{\frac{2n}{n-2}}dz 
 \leq C |\nabla K (x_i)| u_i(x_i)^{-\frac{2}{n-2}}+ Cu_i(x_i)^{-\frac{2(n-1)}{n-2}}.
\end{equation}

\noindent
Plugging this back into (\ref {eq:delta-coarse-pos}) we now have a refined estimate for $\delta_i$:
\begin{equation}
\label {eq:delta-refined-pos}
\delta _i \leq C \left ( u_i(x_i)^{2t_i} +|\nabla K (x_i)| u_i(x_i)^{-\frac{2}{n-2}}+ u_i(x_i)^{-\frac{2(n-1)}{n-2}}\right ).
\end{equation}

To prove (\ref {eq:key-X(K)-pos}) we still need to refine the estimate for $|\nabla K (x_i)|$.

\noindent
In (\ref {eq:gradK-coarse-pos}) we have $$\bigg| \int _{B_{\sigma}} \lambda ^{\frac{2n}{n-2}} u_i^{p_i+1}\frac{\partial K}{\partial z^1} dz \bigg| \leq C \left ( u_i(x_i)^{2 t _i} + \delta _i \right ).$$

\noindent
Again we write out the Taylor expansion
\begin{eqnarray*}
\frac {\partial K}{\partial z^1}(z) & = & \frac {\partial K}{\partial z^1}(x_i) + \sum _{|\beta|=1} \frac{\partial ^{\beta}}{\partial z^{\beta}}\frac {\partial K}{\partial z^1}(x_i)z^{\beta}+\frac{1}{2!}\sum _{|\beta|=2} \frac{\partial ^{\beta}}{\partial z^{\beta}}\frac {\partial K}{\partial z^1}(x_i)z^{\beta}+\cdot \cdot \cdot \\
& & + \frac{1}{(n-3)!}\sum _{|\beta|=n-3} \frac{\partial ^{\beta}}{\partial z^{\beta}}\frac {\partial K}{\partial z^1}(x_i)z^{\beta}+\frac{1}{(n-2)!}\sum _{|\beta|=n-2} \frac{\partial ^{\beta}}{\partial z^{\beta}}\frac {\partial K}{\partial z^1}(\varsigma)z^{\beta}.
\end{eqnarray*}

\noindent
Therefore we have 
\begin{eqnarray*}
& & \int _{B_{\sigma}}\lambda ^{\frac{2n}{n-2}} u_i^{p_i+1} \bigg | \frac{\partial K}{\partial z^1}(x_i) \bigg | dz \\
& \leq & \bigg | \int _{B_{\sigma}} \lambda ^{\frac{2n}{n-2}}u_i^{p_i+1}\frac{\partial K}{\partial z^1} dz \bigg |  + C\sum _{|\beta|=1}^{n-3} \int _{B_{\sigma}} \bigg | \frac{\partial ^{\beta}}{\partial z^{\beta}}\frac {\partial K}{\partial z^1}(x_i)\bigg| |z|^{|\beta|} u_i^{p_i+1}dz\\
& & +C \sum _{|\beta|=n-2} \int _{B_{\sigma}} \bigg | \frac{\partial ^{\beta}}{\partial z^{\beta}}\frac {\partial K}{\partial z^1}(\varsigma)\bigg| |z|^{n-2} u_i^{p_i+1}dz \\
 & \leq &  C \left (u_i(x_i)^{2 t _i} + \delta _i \right ) + C\sum _{|\beta|=1}^{n-3} \int _{B_{\sigma}} \bigg | \frac{\partial ^{\beta}}{\partial z^{\beta}}\frac {\partial K}{\partial z^1}(x_i)\bigg| |z|^{|\beta|} u_i^{p_i+1}dz\\
& & +C \sum _{|\beta|=n-2} \int _{B_{\sigma}} \bigg | \frac{\partial ^{\beta}}{\partial z^{\beta}}\frac {\partial K}{\partial z^1}(\varsigma)\bigg| |z|^{n-2} u_i^{p_i+1}dz. 
\end{eqnarray*}

\noindent
By (\ref {eq:lower-int-u^p+1-pos}) this implies 

\begin{eqnarray*}
\bigg | \frac{\partial K}{\partial z^1}(x_i) \bigg | &  \leq &   C \left ( u_i(x_i)^{2 t _i} + \delta _i \right ) + C\sum _{|\beta|=1}^{n-3} \int _{B_{\sigma}} \bigg | \frac{\partial ^{\beta}}{\partial z^{\beta}}\frac {\partial K}{\partial z^1}(x_i)\bigg| |z|^{|\beta|} u_i^{p_i+1}dz \\
& & +C \sum _{|\beta|=n-2} \int _{B_{\sigma}} \bigg | \frac{\partial ^{\beta}}{\partial z^{\beta}}\frac {\partial K}{\partial z^1}(\varsigma)\bigg| |z|^{n-2} u_i^{p_i+1}dz. 
\end{eqnarray*}

\noindent
By Lemma \ref{lemma:intu^p+1z^k-pos}, (\ref {eq:u^delta-pos}), condition $(**)$, and Young's Inequality, when $1 \leq |\beta| \leq n-3$,
\begin{eqnarray*}
& & \int _{B_{\sigma}} \bigg | \frac{\partial ^{\beta}}{\partial z^{\beta}}\frac {\partial K}{\partial z^1}(x_i)\bigg| |z|^{|\beta|} u_i^{p_i+1}dz  \\
 & \leq & C \int _{B_{\sigma}}  | \nabla K(x_i) |^{\frac{n-1-(|\beta|+1)}{n-2}}   |z|^{|\beta|} u_i^{p_i+1}dz \nonumber \\
& = & C \int _{B_{\sigma}} | \nabla K(x_i) |^{\frac{n-2-|\beta|}{n-2}}\cdot |z|^{|\beta|-\frac{1}{n-2}}\cdot |z|^{\frac{1}{n-2}} u_i^{p_i+1} dz  \\
& \leq & C \int _{B_{\sigma}}\left ( | \nabla K(x_i) |^{\frac{n-2-|\beta|}{n-2}\cdot \frac{n-2}{n-2-|\beta|}}+ |z|^{(|\beta|-\frac{1}{n-2})\cdot \frac{n-2}{|\beta|} } \right )|z|^{\frac{1}{n-2}} u_i^{p_i+1} dz  \\
& = & C \left( \int _{B_{\sigma}}  | \nabla K(x_i) | |z|^{\frac{1}{n-2}}u_i^{p_i+1} dz + \int _{B_{\sigma}} |z|^{n-2-\frac{1}{|\beta|}+\frac{1}{n-2}}u_i^{p_i+1} dz  \right) \\
& \leq &  C |\nabla K(x_i) | u_i(x_i)^{-\frac{2}{(n-2)^2}} + Cu_i(x_i)^{-2+\frac{2}{n-2}\left ( \frac{1}{|\beta|}-\frac{1}{n-2} \right ) }\\
& \leq &  C| \nabla K(x_i) | u_i(x_i)^{-\frac{2}{(n-2)^2}} + Cu_i(x_i)^{-2+\frac{2}{n-2}\frac{n-3}{n-2}},   
\end{eqnarray*}
where the last step holds because $ \frac{1}{|\beta|}-\frac{1}{n-2} \leq  \frac{n-3}{n-2}$.

\noindent
Furthermore, $$\sum _{|\beta|=n-2} \int _{B_{\sigma}} \bigg | \frac{\partial ^{\beta}}{\partial z^{\beta}}\frac {\partial K}{\partial z^1}(\varsigma)\bigg| |z|^{n-2} u_i^{p_i+1}dz \,\,  \leq \,\, C\int _{B_{\sigma}}  |z|^{n-2} u_i^{p_i+1}dz \,\, \leq \,\,  Cu_i(x_i)^{-2}. $$

\noindent
Therefore
\begin{eqnarray*}
& & \bigg | \frac{\partial K}{\partial z^1}(x_i) \bigg | \\
&  \leq &   C \left ( u_i(x_i)^{2 t _i} + \delta _i \right )
+ \left ( C| \nabla K(x_i) | u_i(x_i)^{-\frac{2}{(n-2)^2}} + Cu_i(x_i)^{-2+\frac{2}{n-2}\frac{n-3}{n-2}} \right )  + Cu_i(x_i)^{-2}\\
& \leq & C \delta _i + C u_i(x_i)^{2t _i} +  C| \nabla K(x_i) | u_i(x_i)^{-\frac{2}{(n-2)^2}} + Cu_i(x_i)^{-2+\frac{2}{n-2}\frac{n-3}{n-2}}. 
\end{eqnarray*}

\noindent
The same estimate also holds for $\big | \frac{\partial K}{\partial z^j}(x_i) \big |$, where $j=2,..., n$, so we know

\begin{eqnarray*}
|\nabla K(x_i)  | & \leq &  C \delta _i + C u_i(x_i)^{2t _i} +  C| \nabla K(x_i) | u_i(x_i)^{-\frac{2}{(n-2)^2}} + Cu_i(x_i)^{-2+\frac{2}{n-2}\frac{n-3}{n-2}}   \\
&  \leq &  C \left ( u_i(x_i)^{2t_i} +|\nabla K (x_i)| u_i(x_i)^{-\frac{2}{n-2}}+ u_i(x_i)^{-\frac{2(n-1)}{n-2}} \right ) \\
& & +  C u_i(x_i)^{2t _i} +  C| \nabla K(x_i) | u_i(x_i)^{-\frac{2}{(n-2)^2}} + Cu_i(x_i)^{-2+\frac{2}{n-2}\frac{n-3}{n-2}}    \hspace{.4in} \text{( by (\ref {eq:delta-refined-pos}) )}.
\end{eqnarray*}

\noindent
When $i$ is large enough, all the terms involving  $|\nabla K (x_i)|$ can be absorbed into the left hand side of this inequality, therefore we get a refined estimate 
\begin{equation}
\label {eq:gradK-refined-pos}
 |\nabla K(x_i)|  \leq  C u_i(x_i)^{2t_i} + C u_i(x_i)^{-\frac{2(n-1)}{n-2}}+  Cu_i(x_i)^{-2+\frac{2}{n-2}\frac{n-3}{n-2}}. 
\end{equation}

Finally, we are going to prove (\ref {eq:key-X(K)-pos}).  

\begin{eqnarray*}
& &  u^2_i(x_i)\int _{B_{\sigma}} | X(K)| u_i^{p_i+1}\lambda ^{\frac{2n}{n-2}}dz \\
& = & u^2_i(x_i)\int _{B_{\sigma}} \bigg | r\frac{\partial K}{\partial r}  \bigg| u_i^{p_i+1}\lambda ^{\frac{2n}{n-2}}dz \\
& \leq & C u^2_i(x_i) \left (|\nabla K (x_i)| u_i(x_i)^{-\frac{2}{n-2}}+ u_i(x_i)^{-\frac{2(n-1)}{n-2}} \right ) \hspace{.2in} \text {\big ( by (\ref {eq:int-rdK/dr-u^p+1-pos}) \big )} \\
& \leq &   C u^2_i(x_i) \Bigg (  \left ( u_i(x_i)^{2t_i} + u_i(x_i)^{-\frac{2(n-1)}{n-2}}+ u_i(x_i)^{-2+\frac{2}{n-2}\frac{n-3}{n-2}} \right ) u_i(x_i)^{-\frac{2}{n-2}}  + u_i(x_i)^{-\frac{2(n-1)}{n-2}} \Bigg )\\
& &   \left ( \text {  by  (\ref {eq:gradK-refined-pos}) } \right ) \\
& = & C \left ( u_i(x_i) ^{2+2t_i-\frac{2}{n-2}} + u_i(x_i)^{-\frac{4}{n-2}} + u_i(x_i) ^{-\frac{2}{(n-2)^2}} + u_i(x_i)^{-\frac{2}{n-2}} \right  ).
\end{eqnarray*}

\noindent
By (\ref  {eq: lim-ti-pos}) we know $$\displaystyle \lim _{i \to \infty} \left (2+2t_i-\frac{2}{n-2}\right ) \,\, = \,\, 2 + \frac {2(1-n)}{n} -\frac{2}{n-2} \,\, = \,\, \frac{2}{n}-\frac{2}{n-2} \,\, < \,\, 0,$$ therefore $\displaystyle \lim _{i \to \infty}  u_i(x_i) ^{2+2t_i-\frac{2}{n-2}}=0$.  Then since 
 $$\displaystyle \lim _{i \to \infty} u_i(x_i)^{-\frac{4}{n-2}} =\displaystyle \lim _{i \to \infty} u_i(x_i) ^{-\frac{2}{(n-2)^2}} = \displaystyle \lim _{i \to \infty} u_i(x_i)^{-\frac{2}{n-2}} =0,$$ we have 
$$ \lim _{i \to \infty} u^2_i(x_i)\int _{B_{\sigma}} | X(K)| u_i^{p_i+1}\lambda ^{\frac{2n}{n-2}}dz =0. $$
This proves (\ref {eq:key-X(K)-pos}).  It follows that the limit of the interior term of (\ref {eq:pohozaev-divide-pos}) as $i$ goes to infinity is greater than or equal to $0$.  But this is a contradiction because we have shown that the limit of the boundary term is strictly negative.  Therefore, {\it at least one of the isolated blow-up points must be non-simple.}

\subsubsection{Isolated but Non-simple Blow-up}
\label {subsubsection:non-simple}

Without loss of generality we assume $P_1$ is not a simple blow-up point.Then as a function of $|z|$, $|z|^{\frac{2}{p_i-1}}\bar{u}_i(|z|)$ has a
second critical point at $|z|=r_i$ where $r_i \rightarrow 0$. Let
$y=\frac{z}{r_i}$ and define
$v_i(y)=r_i^{\frac{2}{p_i-1}}u_i(r_iy)$. Then $v_i(y)$ satisfies
\begin{equation}
\label{eq:rescaled-case1-pos}
\Delta _{g^{(i)}} v_i-c(n)R(g^{(i)})v_i+K_iv_i^{p_i}=0 
\end{equation} 
where $g^{(i)}(y) =g_{\alpha \beta}(r_iy)dy^{\alpha}dy^{\beta}$, $R(g^{(i)})(y)=r_i^2R(g)(r_iy)$ and
$K_i(y)=K(r_iy)$.\\ 

\noindent
By this definition $|y|=1$ is the second critical point of $|y|^{\frac{2}{p_i-1}}\bar{v}_i(|y|)$. 
Just as in the scalar-flat case, it can be shown that $0$ is a simple blow-up point for $\{v_i\}$.

By some calculations which are very similar to the proof of Proposition \ref {prop:simpleestimates-pos}, we can prove the following estimates: there exist a constant $C$ independent of
$i$ and a radius $\tilde {r} \leq 1$ such that
\begin{itemize}
\item if \,\, $0
\leq |y| \leq \tilde {r}$, \hspace{.1in} then $$v_i(y) \geq C v_i(0)\left (1+ \frac{K_i(0)}{n(n-2)}v_i(0)^{\frac{4}{n-2}}|y|^2 \right )^{-\frac{n-2}{2}}$$ 
\item  if \,\,  $0 \leq |y| \leq \frac{R}{v_i(0)^{\frac{p_i-1}{2}}}$, \hspace{.1in}
then $$v_i(y) \leq C v_i(0) \left (1+ \frac{K_i(0)}{n(n-2)}
v_i(0)^{p_i-1}|y|^2 \right )^{-\frac{n-2}{2}}$$ 
\item if \,\, $\frac{R}{v_i(0)^{\frac{p_i-1}{2}}} \leq |y| \leq \tilde {r}$,
\hspace{.1in} then \hspace{.1in} $v_i(y) \leq Cv_i(0)^{t _i}|y|^{-l_i}$\\
 where $l_i$, $t_i$ are so chosen that 
$\frac{(2n-1)(n-2)}{2n} < \displaystyle \lim_{i \to \infty} l_i < n-2$,
 and $t _i =1-\frac{(p_i-1)l_i}{2}$.
\end{itemize}

It follows that when $\sigma < \tilde {r}$ and $0 \leq \kappa \leq n-1$, there exists a constant $C$ such that
\begin{equation}
\label {eq:intv^p+1y^k-pos}
\int _{|y|\leq \sigma} |y|^{\kappa}v_i(0)^{p_i+1}dy \leq C v_i(0)^{-\frac{2\kappa}{n-2}+\frac{n-2+\kappa}{2}\delta_i}.
\end{equation}

\noindent
This can be proved by the same calculation as in the proof of Lemma \ref {lemma:intu^p+1z^k-pos}.  Next by an argument that is almost identical to the proof of Proposition \ref {prop:v^delta}, we can show that 

\begin{equation*}
\delta _i \leq C \left ( v_i(0)^{2t _i} + \bigg |\int  _{B_{\sigma}}   \frac {\partial K_i}{\partial r}|y| v_i^{p_i+1}\lambda_i^{\frac{2n}{n-2}}dy  \bigg|\right ) .
\end{equation*} 

\noindent
This gives a preliminary estimate 

$$ \delta _i \leq C \left ( v_i(0)^{-\frac{2}{n-2}}+v_i(0)^{2t _i} \right ),$$and additionally $\displaystyle \lim _{i \to \infty} v_i(0)^{\delta _i}=1$.  Then by the same calculations as those in Section \ref {subsection:prelim-gradK}, we know that for $j=1,2...,n$,

\begin{equation*}
\bigg| \int _{B_{\sigma}} \lambda _i^{\frac{2n}{n-2}}v_i^{p_i+1}\frac{\partial K_i}{\partial
y^j} dy \bigg| \leq C \left ( \delta _i r_i + v_i(0)^{2 t _i} \right ),
\end{equation*}
and we have a preliminary estimate 
$$\bigg |\frac{\partial K_i}{\partial y^j}(0) \bigg | \leq C
\left (r_iv_i(0)^{-\frac{2}{n-2}} + v_i(0)^{2 t _i}
\right ).$$

Choose a point $\bar{y}$ with  $|\bar{y}|=\tilde {r}$.  We have 
$$\Delta _{g^{(i)}} \frac {v_i}{v_i(\bar{y})}-c(n)R(g^{(i)}) \frac {v_i}{v_i(\bar{y})}+ v_i(\bar{y})^{p_i-1} K_i \left (\frac {v_i}{v_i(\bar{y})}\right )^{p_i}=0 .$$

\noindent
On any compact subset $\Omega$ of $\mathbf{R}^n \setminus \{0\}$ which contains $\bar{y}$, since we have a Harnack inequality for $v_i$, $\frac{v_i}{v_i(\bar{y})}$ is uniformly bounded.  Thus because $v_i(\bar{y}) \to 0$ and $g^{(i)}$ converges to the Euclidean metric, $\frac {v_i}{v_i(\bar{y})}$ converges on $\Omega$ in $C^2$-norm to a function $h$ with $\Delta h=0$, where $\Delta$ is the Euclidean Laplacian.  Since $\Omega$ is arbitrary,  $\Delta h=0$ on $\mathbf{R}^n \setminus \{0\}$.  Then because $0$ is a simple blow-up point of $\{v_i\}$ and $|y|^{\frac{2}{p_i-1}}\frac{\bar{v}_i(|y|)}{v_i(\bar{y})}$ has a second critical point at $|y|=1$, we know $h(y)=\frac{1}{2}+\frac{1}{2}|y|^{2-n}$.

Now as in Section \ref {subsection:location} we can prove that $\nabla K(P_1)=\displaystyle \lim _{i \to \infty} \nabla K (x_i)=0$, i.e., $P_1$ is a critical point of $K$.  Recall that the proof is by contradiction: suppose $\nabla K(P_1) \neq 0$, we study the Pohozaev identity (divided by $v_i^2(\bar{y})$) with $X=r\frac{\partial}{\partial r}$ and compare the signs of the limits of both sides.  The key point is to establish the limit
$$\lim _{i \to \infty} v_i^2(0)\bigg | \int _{B_{\sigma}} r\frac{\partial K_i}{\partial r} v_i^{p_i+1} \lambda _i ^{\frac{2n}{n-2}}  dy \bigg |=0. $$

\noindent
In fact, if we have this limit, then by the same argument as in Section \ref {subsection:refined}, it will give a contradiction and rule out Case I completely.

Since $P_1$ is a critical point and $K$ satisfies condition $(**)$, we know 
$$\bigg | \frac{\partial ^{\alpha} K}{\partial z^{\alpha}}(x_i) \bigg | \leq C | \nabla K(x_i) |^{\frac{n-1-|\alpha|}{n-2}} $$ when $  2\leq |\alpha|\leq n-2 .$  Then because $K_i(y)=K(r_iy)$,  
\begin{eqnarray*}
\bigg | \frac{\partial ^{\alpha} K_i}{\partial y^{\alpha}}(0) \bigg | & \leq & C r_i^{\frac{(n-1)(|\alpha| -1)}{n-2}}  | \nabla K_i(0) |^{\frac{n-1-|\alpha|}{n-2}} \\
& < & C r_i  | \nabla K_i(0) |^{\frac{n-1-|\alpha|}{n-2}} \\
& < & C r_i | \nabla K_i(0) |^{\frac{n-2-|\alpha|}{n-3}}, 
\end{eqnarray*}
where the last step uses the fact that $ | \nabla K_i(0) | \to 0$ and $\frac{n-1-|\alpha|}{n-2} > \frac{n-2-|\alpha|}{n-3}$.  Then we can use exactly the same argument as in Section \ref {subsection:refined} to refine the estimates for $\delta _i$ and $|\nabla K_i(0)|$ and thus prove the key limit.  This finishes the proof in Case I.

\subsection{Ruling out Case II}
\label{subsection:case2-pos}
Recall that by defining  $v_i(y)=\sigma _i ^{\frac{2}{p_i-1}}u_i(\sigma _i y)$ and $y=\frac{z}{\sigma _i}$, we have reduced Case II to the situation that  $v_i$ satisfies 
$$
\Delta _{g^{(i)}} v_i -c(n)R \left (g^{(i)}\right )v_i+ K_iv_i^{p_i}=0
$$
where
$g^{(i)}(y)=g_{\alpha \beta}(\sigma _i y)dy^{\alpha}dy^{\beta}$, $R(g^{(i)})(y)=\sigma _i ^2 R(g)(\sigma _iy)$ and $K_i(y)=K(\sigma _iy)$, and $0$ is an isolated blow-up point of $\{v_i\}$.  

If $0$ is not a simple blow-up point, then we can do another rescaling and repeat the previous argument in Section \ref {subsubsection:non-simple}, with $r_i$ replaced by $r_i\sigma _i$, to get a contradiction.  Therefore $0$ must be a simple blow-up point for $\{v_i\}$.  Then we can still repeat the argument in Section \ref {subsubsection:non-simple}, with $r_i$ replaced by $\sigma _i$. The only difference is in 
the expression of $h=\displaystyle \lim _{i \rightarrow \infty}
\frac{v_i(y)}{v_i(\bar{y})}$.  As in the scalar-flat case, because here $|y|^{\frac{2}{p_i-1}} \bar{v}_i(|y|)$ doesn't have a second critical point at $|y|=1$, we have a different expression of $h$: near $0$,
$$h(y) = c_1|y|^{2-n}+  A + O(|y|) $$ where $A$ is a positive constant.  This positive ``mass'' term $A >0 $ guarantees that the limit of the boundary term of the Pohozaev identity is still negative, i.e., $$\lim _{i
\rightarrow \infty} \frac{1}{v^2_i(\bar{y})}\int _{\partial
B_{\sigma}} T_i(X, \nu _i) d \Sigma _i <0.$$  The other parts of the proof remain the same.  Therefore Case II can also be ruled out.

This completes the proof of Theorem \ref {thm:main-pos}.

\bibliographystyle{plain}
    \bibliography{thesis}

\end{document}